\documentclass[11pt]{article}
\usepackage{cmap}
\usepackage[T1]{fontenc}
\usepackage{lmodern}
\usepackage{amsmath,amssymb,amsthm}
\usepackage{mathtools,mathrsfs}
\usepackage{array}
\usepackage{booktabs}
\usepackage{graphicx}
\usepackage[margin=1.15in]{geometry}
\usepackage{microtype}
\usepackage{enumitem}
\usepackage[section]{placeins}
\usepackage[font=small,labelfont=bf,labelsep=period]{caption}
\usepackage[dvipsnames]{xcolor}
\usepackage[
  backend=biber,
  style=alphabetic,
  sorting=nyt,
  maxbibnames=99,
  maxalphanames=4,
  natbib=true
]{biblatex}
\appto{\bibsetup}{\emergencystretch=1em}
\usepackage{aliascnt}
\usepackage{hyperref}
\usepackage[nameinlink,noabbrev,capitalise]{cleveref}
\graphicspath{{figures/}}
\numberwithin{equation}{section}

\setlist[itemize]{leftmargin=2em,itemsep=0.3em,topsep=0.4em}
\setlist[enumerate]{leftmargin=2em,itemsep=0.3em,topsep=0.4em}

\newtheorem{theorem}{Theorem}[section]
\newaliascnt{proposition}{theorem}
\newtheorem{proposition}[proposition]{Proposition}
\aliascntresetthe{proposition}
\newaliascnt{lemma}{theorem}
\newtheorem{lemma}[lemma]{Lemma}
\aliascntresetthe{lemma}
\newaliascnt{corollary}{theorem}
\newtheorem{corollary}[corollary]{Corollary}
\aliascntresetthe{corollary}
\newaliascnt{conjecture}{theorem}

\aliascntresetthe{conjecture}
\theoremstyle{definition}
\newaliascnt{definition}{theorem}

\aliascntresetthe{definition}
\newaliascnt{assumption}{theorem}

\aliascntresetthe{assumption}
\newaliascnt{example}{theorem}
\newtheorem{example}[example]{Example}
\aliascntresetthe{example}
\theoremstyle{remark}
\newaliascnt{remark}{theorem}
\newtheorem{remark}[remark]{Remark}
\aliascntresetthe{remark}
\crefname{theorem}{\mbox{Theorem}}{\mbox{Theorems}}
\Crefname{theorem}{\mbox{Theorem}}{\mbox{Theorems}}
\crefname{proposition}{\mbox{Proposition}}{\mbox{Propositions}}
\Crefname{proposition}{\mbox{Proposition}}{\mbox{Propositions}}
\crefname{lemma}{\mbox{Lemma}}{\mbox{Lemmas}}
\Crefname{lemma}{\mbox{Lemma}}{\mbox{Lemmas}}
\crefname{corollary}{\mbox{Corollary}}{\mbox{Corollaries}}
\Crefname{corollary}{\mbox{Corollary}}{\mbox{Corollaries}}
\crefname{conjecture}{\mbox{Conjecture}}{\mbox{Conjectures}}
\Crefname{conjecture}{\mbox{Conjecture}}{\mbox{Conjectures}}
\crefname{definition}{\mbox{Definition}}{\mbox{Definitions}}
\Crefname{definition}{\mbox{Definition}}{\mbox{Definitions}}
\crefname{assumption}{\mbox{Assumption}}{\mbox{Assumptions}}
\Crefname{assumption}{\mbox{Assumption}}{\mbox{Assumptions}}
\crefname{example}{\mbox{Example}}{\mbox{Examples}}
\Crefname{example}{\mbox{Example}}{\mbox{Examples}}
\crefname{remark}{\mbox{Remark}}{\mbox{Remarks}}
\Crefname{remark}{\mbox{Remark}}{\mbox{Remarks}}
\crefname{section}{\mbox{Section}}{\mbox{Sections}}
\Crefname{section}{\mbox{Section}}{\mbox{Sections}}
\crefname{subsection}{\mbox{Section}}{\mbox{Sections}}
\Crefname{subsection}{\mbox{Section}}{\mbox{Sections}}
\crefname{appendix}{\mbox{Appendix}}{\mbox{Appendices}}
\Crefname{appendix}{\mbox{Appendix}}{\mbox{Appendices}}
\crefname{subappendix}{\mbox{Appendix}}{\mbox{Appendices}}
\Crefname{subappendix}{\mbox{Appendix}}{\mbox{Appendices}}
\crefname{figure}{\mbox{Figure}}{\mbox{Figures}}
\Crefname{figure}{\mbox{Figure}}{\mbox{Figures}}
\crefname{table}{\mbox{Table}}{\mbox{Tables}}
\Crefname{table}{\mbox{Table}}{\mbox{Tables}}

\newcommand{\R}{\mathbb R}
\newcommand{\Nzero}{\mathbb N_0}
\newcommand{\supp}{\operatorname{supp}}
\newcommand{\aff}{\operatorname{aff}}
\newcommand{\Span}{\operatorname{span}}
\newcommand{\Cov}{\operatorname{Cov}}
\newcommand{\diag}{\operatorname{diag}}

\newcommand{\cof}{\operatorname{cof}}
\newcommand{\inter}{\operatorname{int}}
\newcommand{\dd}{\mathrm d}
\newcommand{\ee}{\mathrm e}
\newcommand{\KL}{D_{\mathrm{KL}}}

\title{On Natural Exponential Families with\\
Log-Affine Generalized Variance}
\author{%
Zeyu Cao\thanks{Independent Researcher, Long Island City, NY 11101, USA.
Email: \texttt{josephcao891011@gmail.com}.}
\and
Shaosai Huang\thanks{Kspectra Research Inc., Toronto, ON, M2N 0G3, Canada.
Email: \texttt{arthur.foxie.huang@kspectra.ai}.}}
\date{November 3, 2026}
\hypersetup{
  colorlinks=true,
  linkcolor=MidnightBlue,
  citecolor=OliveGreen,
  urlcolor=BrickRed,
  pdftitle={On Natural Exponential Families with Log-Affine Generalized Variance},
  pdfauthor={Zeyu Cao and Shaosai Huang},
  pdfsubject={Regularity, Gaussian rigidity and boundary structure of natural exponential families with log-affine generalized variance},
  pdfkeywords={constant regression, convex support, generalized variance, infinite divisibility, Monge–Ampère equation, natural exponential family, Poisson–Gaussian family, projectively induced Kähler metric, semi-flat Calabi–Yau metric},
  pdflang={en-US}
}

\begin{document}
\maketitle

\begin{abstract}
We study natural exponential families whose generalized variance, the
determinant of the covariance matrix, is log-affine in the natural
parameter.  Without any assumption on the natural parameter domain, we show
that every such family is regular and that it is Gaussian when its
generalized variance is constant.  Consequently, a sample of size $n>d$ in $\R^d$ comes from a
nondegenerate normal law exactly when the conditional expectation of the
sample generalized variance given the sample mean is a positive constant.
Every non-Gaussian family puts mass on a facet of its convex
support, and, as for a Poisson count, its height above that facet takes only
integer multiples of a fixed step.  More generally, its displacement away
from any exposed face that carries mass has nonnegative integer coordinates
in a suitable basis.  A local inequality between third and
fourth cumulants singles out the Poisson--Gaussian families and yields a
complete classification of the infinitely divisible ones.  Geometrically, a
semi-flat K\"ahler tube, whose potential has constant Hessian determinant,
is flat if some multiple of its metric is projectively induced near a
point.
\end{abstract}

\noindent\textit{2020 Mathematics Subject Classification.} Primary 62E10,
60E05; Secondary 60E07, 62H05, 35J96, 52A20, 53C55.

\noindent\textit{Keywords.} Constant regression; convex support; generalized
variance; infinite divisibility; Monge--Amp\`ere equation; natural
exponential family; Poisson--Gaussian family; projectively induced K\"ahler
metric; semi-flat Calabi--Yau metric.

\section{Introduction}
\label{sec:intro}

The generalized variance of a random vector is the determinant of its
covariance matrix.  For a natural exponential family it is the Hessian
determinant of the cumulant function, and its square root is proportional
to the density of the Jeffreys prior.  This paper studies the natural
exponential families whose generalized variance is the exponential of an
affine function of the natural parameter, that is, whose Jeffreys prior is
log-affine.

To fix notation, let $\mu$ be a positive Radon measure on $\R^d$ whose
support is not contained in an affine hyperplane; we say that such a measure
has \emph{full affine support}.  Vectors are column vectors: a tuple such
as $(s,t)$ stands for the column obtained by stacking its
entries, ${}^{\mathsf T}$ denotes the transpose, and
$\langle u,v\rangle=u^{\mathsf T}v$.  Its Laplace transform and its natural
parameter domain are
\[
 Z(\theta)=L_\mu(\theta)=\int_{\R^d}\ee^{\langle\theta,x\rangle}\,\mu(\dd x),
 \qquad
 D=\{\theta\in\R^d:Z(\theta)<\infty\}.
\]
Assume that the interior $\Theta$ of the convex set $D$ is nonempty and
let $K=\log Z$ on $\Theta$ be the cumulant function.  The natural
exponential family (NEF) generated by $\mu$ consists of the laws
\[
 P_\theta(\dd x)=\exp\{\langle\theta,x\rangle-K(\theta)\}\,\mu(\dd x),
 \qquad\theta\in\Theta.
\]
The law $P_\theta$ has mean $\nabla K(\theta)$ and covariance matrix
$\nabla^2K(\theta)$, which is also the Fisher information of the natural
parameter; its determinant $\det\nabla^2K(\theta)$ is the generalized
variance.
Write $M=\nabla K(\Theta)$ for the mean domain and $Q$ for the closed convex
hull of $\supp\mu$, the closed convex support.  The families we study are
those for which
\begin{equation}\label{eq:E}
 \det\nabla^2K(\theta)=C\ee^{\langle b,\theta\rangle}
 \qquad(\theta\in\Theta),
 \qquad C>0,\quad b\in\R^d.
\end{equation}
We then say that the family has \emph{exponential-affine generalized
variance} and call $\mu$ a \emph{solution} of \eqref{eq:E}.  The case $b=0$
is exactly the case of a uniform Jeffreys prior.

Equation \eqref{eq:E} is invariant under invertible affine changes
$x\mapsto Ax+c$ of the statistic (which replace $b$ by $Ab$), under
exponential tilts $\mu\mapsto\ee^{\langle\theta_0,x\rangle}\mu$ and under
multiplication of $\mu$ by positive constants.  Its basic solutions are the
Gaussian families, for which $b=0$, and the Poisson--Gaussian families: in
suitable affine coordinates, writing $\theta=(s,t)\in\R^r\times\R^k$,
\begin{equation}\label{eq:PG}
 K(s,t)=\tfrac12s^{\mathsf T}\Sigma s+\langle c,(s,t)\rangle
 +\sum_{j=1}^k\lambda_j(\ee^{t_j}-1),
 \qquad \Sigma>0,\ \lambda_j>0,
\end{equation}
which satisfies \eqref{eq:E} with $b=(0,\dots,0,1,\dots,1)$.  These are the
families generated by the product of a nondegenerate Gaussian law on $\R^r$
and $k$ independent Poisson laws, or by a positive multiple of it; a
multiple adds a constant to $K$ and changes neither the laws $P_\theta$ nor
\eqref{eq:E}.  We call a
solution \emph{Gaussian} if $\mu$ is a positive multiple of a
nondegenerate normal distribution, and \emph{Poisson--Gaussian} if, after an
invertible affine change of the statistic, $K$ takes the form \eqref{eq:PG}
up to an additive constant, that is, if $\mu$ is a positive multiple of an
affine image of such a product law; the Gaussian solutions are the case
$k=0$.

Whether \eqref{eq:E} has other solutions is a question from the theory of
conjugate priors.  For NEFs on $\R^d$, the standard conjugate families of
the natural and of the mean parametrization coincide exactly when
$\det\nabla^2K=\exp\{\langle a,\theta\rangle+\beta K(\theta)+\gamma\}$
(see \cite{GutierrezPenaSmith1995} and \cite[p.~1831]{Casalis1996}), and
Casalis left the determination of all such families open; \eqref{eq:E} is
the case $\beta=0$.  The Poisson--Gaussian families are the only known
solutions, and it is expected that there are no others.

The earlier partial answers rest on a link between generalized variance and
Laplace transforms.  Kokonendji and Seshadri~\cite{KokonendjiSeshadri1996} showed
that $\det\nabla^2K$, multiplied by a power of $Z$, is again a Laplace
transform (\cref{lem:KS} below), and proved that constant generalized
variance forces a Gaussian family when $D=\R^d$
\cite[Theorem~3.5, pp.~1825--1826]{KokonendjiSeshadri1996}.  For infinitely
divisible (ID) measures, Hassairi~\cite[Theorem~2.1]{Hassairi1999} showed
that $\det\nabla^2K$ itself is a Laplace transform.  Kokonendji and
Masmoudi~\cite[Theorem~3]{KokonendjiMasmoudi2006} stated that the ID
solutions of \eqref{eq:E} are exactly the Poisson--Gaussian families when
$D=\R^d$, a hypothesis they regarded as crucial; their proof relies on a
lemma that is false as stated (\cref{rem:KM-lemma}).  They described
the case $b=0$ on a proper domain as a challenging problem
\cite[p.~376]{KokonendjiMasmoudi2006}.  Earlier,
Kokonendji~\cite[Theorem~3.2]{Kokonendji2003} had announced the
Poisson--Gaussian classification of \eqref{eq:E} in full generality, but
the printed proof is incomplete (\cref{rem:2003}).  Related
characterizations of other families by their generalized variance are given
in
\cite{GhribiMasmoudi2010,KokonendjiMasmoudi2013,KokonendjiMoypemnaSembonaNisa2020}.

The differential equation alone does not decide the question.  As a
partial differential equation, \eqref{eq:E} has many convex solutions that
are not cumulant functions (\cref{rem:not-cumulant}).  Equation
\eqref{eq:E} and its Legendre dual are the affine K\"ahler--Ricci-flat
equations of Li and Xu \cite[(1.3)--(1.4)]{LiXu2009}.  For entire
solutions, that is, solutions on all of $\R^d$, of these equations and of
coupled variants, see \cite{LiXu2009,Wang2026,SongWuXu2026}; on general
domains, Li and Xu proved rigidity under a completeness assumption on an
associated metric \cite{LiXu2009b}.  In the NEF problem, $\Theta$ may be a
proper subset of $\R^d$ and no completeness is assumed; what replaces these
hypotheses is the positivity of $\mu$.

We settle the case $b=0$.  For $b\ne0$ the classification remains open, but
we show that positivity of $\mu$ constrains every solution at two scales.  Globally,
the determinant transform of \cref{sec:transform} converts
\eqref{eq:E} into an equality of positive measures and hence into
arithmetic information about the boundary of the convex support.  Locally,
differentiating \eqref{eq:E} twice shows that a tensor built from the third
and fourth cumulants has zero trace, so a one-sided sign condition forces
that tensor to vanish.  We call these the boundary-quantization and
cumulant-rigidity mechanisms.  The results below follow the order of the
sections that prove them.

\medskip\noindent\textbf{\textit{Regularity and the Constant Case.}}
Recall that an NEF is \emph{regular} if its natural parameter domain $D$ is
open, and \emph{steep} if $|\nabla K(\theta_n)|\to\infty$ along every
sequence in $\Theta$ that converges to a boundary point of $\Theta$
\cite[Chapter~8]{BarndorffNielsen1978}.  The earlier results on
\eqref{eq:E} recalled above assume $D=\R^d$.  Our first result shows that
regularity and steepness are automatic.

\begin{theorem}[Regularity]\label{thm:regularity-intro}
If $\mu$ satisfies \eqref{eq:E}, then $Z(\theta)=+\infty$ for every
$\theta\in\partial\Theta$.  Hence the family is regular: $D=\Theta$ is
open, $K(\theta)\to+\infty$ as $\theta$ tends to any finite boundary point
of $\Theta$, the family is steep, and $M=\inter Q$.
\end{theorem}

At a hypothetical boundary point where $Z$ is finite, the logarithm of the
tilted characteristic function would satisfy an analytic elliptic equation,
which produces exponential moments.  Regularity makes the graph of $K$
closed.  For $b=0$, the closed-graph form of the J\"orgens--Calabi--Pogorelov
theorem due to Cheng and Yau~\cite{ChengYau1986} then gives the following
answer to the question of Kokonendji and Masmoudi, with no assumption on the
natural parameter domain, on steepness or on infinite divisibility.

\begin{theorem}[Constant Generalized Variance]\label{thm:constant-intro}
If $\det\nabla^2K$ is constant on $\Theta$, then $\mu$ is Gaussian; in
particular $\Theta=\R^d$.
\end{theorem}

In terms of the Jeffreys prior, \cref{thm:constant-intro} reads as
follows.

\begin{corollary}[Jeffreys Prior]\label{cor:jeffreys-intro}
A natural exponential family whose Jeffreys prior is uniform in the natural
parameter is Gaussian.
\end{corollary}

\Cref{thm:constant-intro} also gives a characterization of the
normal law that does not involve exponential families, by a regression
property in the sense of Lukacs~\cite{Lukacs1979}: a statistic $T$ has
\emph{constant regression} on a statistic $U$ if
$\mathbb E[T\mid U]=\mathbb E[T]$ almost surely.  Let $X_1,\dots,X_n$ be
independent, identically distributed random vectors in $\R^d$, with sample
mean and scatter matrix
\[
 \bar X=\frac1n\sum_{j=1}^nX_j,
 \qquad
 W=\sum_{j=1}^n(X_j-\bar X)(X_j-\bar X)^{\mathsf T}.
\]
Up to the factor $(n-1)^{-d}$, $\det W$ is the generalized variance of the
sample; for $n=d+1$ it is $(d!)^2/(d+1)$ times the squared volume of the
simplex with vertices $X_1,\dots,X_{d+1}$.  For a normal sample, $W$ is
independent of $\bar X$, so $\det W$ has constant regression on $\bar X$.
The converse holds as soon as $\det W$ is not almost surely zero, and it
needs no moment assumption: the conditional expectation of the nonnegative
variable $\det W$ is always defined, with values in $[0,\infty]$.

\begin{corollary}[Regression Characterization of the Normal Law]
\label{cor:regression-intro}
Let $n\ge d+1$.  The conditional expectation
$\mathbb E[\det W\mid\bar X]$ is almost surely equal to a positive finite
constant if and only if the $X_j$ have a nondegenerate normal distribution.
\end{corollary}

For $n=d+1$ and under the assumption that the Laplace transform of the law
is finite on all of $\R^d$, this is the theorem of Kokonendji and Seshadri
recalled above, which they state in an equivalent form
\cite[Theorem~3.5]{KokonendjiSeshadri1996}.  The classical theorem of Lukacs
and Laha \cite[Theorem~6.2.9]{LukacsLaha1964}, quoted in
\cite[Theorem~2.2]{Lukacs1979}, assumes a finite covariance matrix and
constant regression on $\bar X$ of the matrix $W$ itself, or of another
quadratic matrix statistic.  \Cref{cor:regression-intro} requires
it of the single scalar $\det W$ instead; for $d=1$, where $\det W=W$, it
is contained in that theorem.
The proof uses that the argument for \cref{thm:regularity-intro}
needs only the measure identity of \cref{sec:transform}, in the
convolved form \eqref{eq:conv-identity}: a probability law that satisfies
it has exponential moments, whether or not it is known to generate an
exponential family (\cref{prop:exp-moments}).

\medskip\noindent\textbf{\textit{Global Boundary Rigidity.}}
For $b\neq0$ the Poisson factors in \eqref{eq:PG} put mass on facets of $Q$
and are integer-valued in the normal direction.  Every non-Gaussian
solution shares this structure.  To state it, recall that an \emph{exposed
face} of the closed convex set $Q$ is a set of the form $Q\cap\{h=0\}$,
where $h$ is an affine function with $h\ge0$ on $Q$; it is \emph{proper} if
it is nonempty and different from $Q$, a \emph{facet} if it has dimension
$d-1$, and an \emph{exposed point} if it is a single point
\cite[Section~18]{Rockafellar1970}.  We write $\aff A$ for the affine hull
of a set $A$ and $f_\#\mu$ for the image of $\mu$ under a map $f$, call a
point of positive mass an \emph{atom}, say that $\mu$ \emph{charges} a set
of positive measure, and put $\Nzero=\{0,1,2,\dots\}$.  For a proper exposed
face $F$ of $Q$ with
$\mu(F)>0$, translate a point of $\supp(\mu|_F)$ to the origin and put
$V=\Span(\supp\mu|_F)$, $N=\R^d/V$, $k=\dim N$, and let
$\pi:\R^d\to N$ be the quotient map, which records the directions normal to
the face.

\begin{theorem}[Boundary Structure]\label{thm:boundary-intro}
Let $\mu$ satisfy \eqref{eq:E} and be non-Gaussian; equivalently,
$b\ne0$.  Then the following hold.
\begin{enumerate}
 \item[\textup{(a)}] $Q$ has a proper exposed facet $F$ with $\mu(F)>0$ and
 $\aff\supp(\mu|_F)=\aff F$.  If
 \begin{equation*}
  F=Q\cap\{h=0\},\quad\text{where } h(x)=\langle a,x\rangle-\alpha\ge0
  \ \text{on } Q,
 \end{equation*}
 then $g=\langle a,b\rangle>0$,
 the measure $h_\#\mu$ is supported in $g\Nzero$, and it has atoms at $0$
 and $g$.
 \item[\textup{(b)}] For every proper exposed face $F$ with $\mu(F)>0$, there
 are atoms $g_1,\dots,g_k$ of $\pi_\#\mu$ forming a basis of $N$ such that
 \[
  \pi(b)=g_1+\dots+g_k,
  \qquad
  \supp(\pi_\#\mu)\subset\Nzero g_1+\dots+\Nzero g_k.
 \]
 \item[\textup{(c)}] If $\mu$ charges an exposed point $v$ of $Q$, there are
 linearly independent $g_1,\dots,g_d$, whose sum is $b$, such that
 $\mu$ charges $v,v+g_1,\dots,v+g_d$ and
 \begin{equation*}
  \supp\mu\subset v+\Nzero g_1+\dots+\Nzero g_d.
 \end{equation*}
 After an affine change
 of coordinates, $\mu$ lives on $\Nzero^d$ and the family consists of
 $d$-variate power-series distributions.
\end{enumerate}
By \textup{(a)}, every solution that is absolutely continuous, or whose
convex support is $\R^d$, is Gaussian.  Moreover, a non-Gaussian solution
has unbounded support.
\end{theorem}

Limits of $P_\theta$ at exposed faces were studied by Csisz\'ar and
Mat\'u\v s~\cite{CsiszarMatus2005}, Hassairi and
Masmoudi~\cite{HassairiMasmoudi2005} and Masmoudi~\cite{Masmoudi2006}.  The
existence of a charged facet uses the Legendre dual of $K$ and a
Liouville-type nonexistence theorem for the dual equation
(\cref{prop:dual-nonexistence}), which was stated by Li and
Xu~\cite{LiXu2009} and is proved here by adapting an argument of
Wang~\cite{Wang2026}; the lattice structure then follows from the support
identity of \cref{sec:transform}.

\medskip\noindent\textbf{\textit{Local Cumulant Rigidity.}}
The second mechanism is local in the parameter.  For $X\sim P_\theta$ and
$u,v\in\R^d$, the joint cumulants of $X$ and of its projections
$Y_u=\langle u,X\rangle$ and $Y_v=\langle v,X\rangle$ are derivatives of $K$
at $\theta$; we write $\operatorname{cum}$ for joint cumulants.  In
dimension one the inequality below says that the squared skewness of
$P_\theta$ is at most its excess kurtosis, which holds for every infinitely
divisible law~\cite{RohatgiSzekely1989}.  Under \eqref{eq:E} it is rigid.

\begin{theorem}[Local Fourth-Cumulant Rigidity]\label{thm:cum-intro}
Let $\mu$ satisfy \eqref{eq:E}.  It is Poisson--Gaussian if and only if there
is a nonempty open $U\subset\Theta$ such that, for $X\sim P_\theta$,
$Y_u=\langle u,X\rangle$ and $Y_v=\langle v,X\rangle$,
\[
 \operatorname{cum}(Y_u,Y_v,Y_u,Y_v)
 \ge c^{\mathsf T}\Cov(X)^{-1}c,
 \qquad c=\operatorname{cum}(Y_u,Y_v,X),
\]
for every $\theta\in U$ and $u,v\in\R^d$.  For every other solution, the set
of $\theta$ for which some $u,v$ give the strict reverse inequality is open
and dense in $\Theta$.
\end{theorem}

Every infinitely divisible family satisfies the inequality of
\cref{thm:cum-intro} on all of $\Theta$, so the theorem contains the
following classification, which needs no assumption on the natural parameter
domain.

\begin{corollary}[Infinitely Divisible Families]\label{cor:id-intro}
If $\mu$ satisfies \eqref{eq:E} and the laws $P_\theta$ are infinitely
divisible, then there is an invertible affine map under which every
$P_\theta$ becomes the product of a nondegenerate Gaussian law on $\R^r$
and $d-r$ independent Poisson laws, for some $0\le r\le d$.
\end{corollary}

This is the classification stated in
\cite[Theorem~3]{KokonendjiMasmoudi2006} for $D=\R^d$.  The criterion rests
on the characterization of the Poisson--Gaussian families as the NEFs whose
variance function, the covariance matrix expressed in terms of the mean, is
affine; this is due to Letac~\cite{Letac1989}.  Furuhata and
Kurose~\cite{FuruhataKurose2013} classified the corresponding Hessian
manifolds, those of constant Hessian sectional curvature zero, and Tojo and
Yoshino~\cite{TojoYoshino2026} treat one-dimensional exponential families.
The trace identity behind \cref{thm:cum-intro} is the standard
expression of Ricci curvature as a trace of the curvature tensor.  What is
new is that under \eqref{eq:E} the inequality is needed only for pairs of
directions and only on some nonempty open set, the open-dense failure
alternative, and the resulting complete proof of
\cref{cor:id-intro}.

The classification also constrains every other solution.  For $r>0$, the
$r$th convolution power of $\mu$ exists if $Z^r$ is again the Laplace
transform of a positive measure.  The set of such $r$, the J{\o}rgensen set
of $\mu$, is all of $(0,\infty)$ exactly when the laws $P_\theta$ are
infinitely divisible, and it is a closed additive semigroup.  Hence a
solution that is not Poisson--Gaussian has no convolution powers of small
order (\cref{sec:cumulant}).

\medskip\noindent\textbf{\textit{Further Global Rigidity in Dimension Two.}}
In the plane, a non-Gaussian solution whose convex support contains no
whole line is, after an affine change of coordinates, a bivariate
power-series family on the lattice orthant $\Nzero^2$, as for two
independent Poisson laws; otherwise its convex support is invariant under
translations along a line, as for the product of a Poisson and a Gaussian
law (\cref{thm:rank-two-intro}).  Neither case is classified here, but a
reduction to positive axisymmetric harmonic functions on $\R^3$ gives
rigidity in three regular regimes (\cref{thm:harmonic-intro}).

\medskip\noindent\textbf{\textit{Geometric Consequences.}}
Equation \eqref{eq:E} expresses the Ricci-flatness of the K\"ahler metric
$\omega_K=i\partial\bar\partial K(\operatorname{Re}z)$ on the tube
$\Theta+i\R^d$ \cite{ZhangKhan2020}.  A K\"ahler metric is
\emph{projectively induced} if it is the pull-back of the Fubini--Study
metric under a holomorphic immersion, called a K\"ahler immersion, into a
finite- or infinite-dimensional complex projective space.  When
$K=\log L_\mu$, the metric $\omega_K$ is projectively induced by the map
$z\mapsto[x\mapsto\ee^{\langle z,x\rangle/2}]$ into
$\mathbb P(L^2(\mu))$, which we call the coherent-state map.
\Cref{sec:kahler} proves a converse and two consequences.  First, for
a non-Poisson--Gaussian solution the rescaled metrics $r\,\omega_K$ are not
projectively induced for small $r>0$, even locally.  Second, let $K$ be any
smooth strictly convex function with $\det\nabla^2K$ constant.  Then the
fibres $\{\theta\}+i\R^d$ of the tube are flat and special Lagrangian, so
that $\omega_K$ is a semi-flat Calabi--Yau
metric~\cite{LeungYauZaslow2000}; if some multiple $c\,\omega_K$ is
projectively induced, even near a single point, then
$\omega_K$ is flat (\cref{thm:semiflat}).  This verifies a further
case of the conjecture of Loi, Salis and Zuddas that Ricci-flat projectively
induced K\"ahler metrics are flat
\cite[Conjecture~1]{LoiSalisZuddas2018}.

\medskip\noindent\textbf{\textit{Scope and Organization.}}
For $b=0$ the results above are complete: every solution of \eqref{eq:E}
is Gaussian.  For $b\ne0$ we prove the boundary structure, the cumulant
criterion with the classification of the infinitely divisible solutions,
and the results in dimension two; whether every solution with $b\ne0$ is
Poisson--Gaussian remains open (\cref{rem:2003}).
Besides the results cited above, the proofs use Morrey's analyticity
theorem~\cite{Morrey1958} and the weak continuity of Monge--Amp\`ere
measures~\cite{Gutierrez2016}; the remaining arguments are elementary.  \Crefrange{sec:transform}{sec:cumulant} hold in every
dimension: \cref{sec:transform} derives the determinant transform and
its support identity, \cref{sec:regularity} proves regularity and
treats the constant case, \cref{sec:boundary} develops the global
boundary mechanism, and \cref{sec:cumulant} the local cumulant
mechanism.  \Cref{sec:kahler} gives the geometric applications, and
\cref{sec:remarks} compares our results with two earlier arguments.
\Cref{app:cumulant} proves the last assertion of
\cref{lem:cum-master}, gives a direct proof of the classification step
in \cref{thm:cum}, and contains an example showing that \eqref{eq:E}
is needed in that theorem.  \Cref{app:dimension-two} states and proves
the results in dimension two.

\section{The Determinant Transform}
\label{sec:transform}

The determinant transform is the algebraic engine of the global argument:
it replaces \eqref{eq:E} by an equality of positive measures and therefore
by an identity of supports.

Replacing $\mu$ by $P_{\theta_0}$ and $\theta$ by $\theta-\theta_0$
preserves \eqref{eq:E}, with $C$ replaced by $C\ee^{\langle
b,\theta_0\rangle}$.  Hence we may assume, when convenient, that
$0\in\Theta$ and $Z(0)=1$.  For $x\in\R^d$ write
$\widetilde x=\bigl(\begin{smallmatrix}1\\ x\end{smallmatrix}\bigr)\in\R^{d+1}$, the column vector with first entry $1$
followed by the entries of $x$.  For points $x_0,\dots,x_d\in\R^d$, write
$\det(\widetilde x_0,\ldots,\widetilde x_d)$ for the determinant of the
$(d+1)\times(d+1)$ matrix with columns $\widetilde x_0,\dots,\widetilde x_d$.
Subtracting the first column from the others shows that it equals
$\det(x_1-x_0,\dots,x_d-x_0)$, that is, $\pm d!$ times the volume of the
simplex with vertices $x_0,\dots,x_d$.  The \emph{determinant transform}
of $\mu$ is the positive measure $\nu_\mu$ on $\R^d$ given, for Borel sets
$E\subset\R^d$, by
\begin{equation}\label{eq:nu-def}
 \nu_\mu(E)=\frac1{(d+1)!}\int_{(\R^d)^{d+1}}
 \det(\widetilde x_0,\ldots,\widetilde x_d)^2\,
 \mathbf 1_E(x_0+\dots+x_d)\,\mu(\dd x_0)\cdots\mu(\dd x_d).
\end{equation}
The integral is over $d+1$ points $x_0,\dots,x_d$ of $\R^d$, each
integrated against $\mu$.  For a probability measure $P$ and independent
random vectors $X_0,\dots,X_d$ with law $P$, $\nu_P(E)$ is therefore
$1/(d+1)!$ times the expectation of
$\det(\widetilde X_0,\ldots,\widetilde X_d)^2$ on the event
$\{X_0+\dots+X_d\in E\}$.
Kokonendji and Seshadri computed its Laplace transform.

\begin{lemma}[Kokonendji--Seshadri]\label{lem:KS}
For $\theta\in\Theta$,
\[
 \int\ee^{\langle\theta,s\rangle}\nu_\mu(\dd s)
 =Z(\theta)^{d+1}\det\nabla^2K(\theta).
\]
Consequently, \eqref{eq:E} holds
if and only if
\begin{equation}\label{eq:measure-identity}
 \nu_\mu=C\,\delta_b*\mu^{*(d+1)}.
\end{equation}
\end{lemma}

This is \cite[Theorem~2.2]{KokonendjiSeshadri1996}; see also
\cite[Proposition~2.4]{Kokonendji2003}.  We recall the short proof.

\begin{proof}
Put
\[
 M(\theta)=\int\ee^{\langle\theta,x\rangle}
 \widetilde x\widetilde x^{\mathsf T}\mu(\dd x).
\]
The continuous Cauchy--Binet formula (Andr\'eief's identity) identifies the
Laplace transform of $\nu_\mu$ with $\det M(\theta)$.  Since
\[
 Z^{-1}M=\begin{pmatrix}1&m^{\mathsf T}\\ m&\nabla^2K+mm^{\mathsf T}
 \end{pmatrix},\qquad m=\nabla K(\theta),
\]
the Schur complement gives $\det M=Z^{d+1}\det\nabla^2K$.  Under
\eqref{eq:E} the right side is the Laplace transform of
$C\delta_b*\mu^{*(d+1)}$, and uniqueness of Laplace transforms on the open
set $\Theta$ gives \eqref{eq:measure-identity}.  The converse follows by
taking transforms.
\end{proof}

The measure identity \eqref{eq:measure-identity} survives exponential
tilting, even at a boundary point of $\Theta$ where $Z$ is finite; the
regularity proof starts from this observation.

\begin{lemma}[Stability under Exponential Tilting]\label{lem:tilt}
Suppose \eqref{eq:measure-identity} holds and $Z(\theta_*)<\infty$ for
some $\theta_*\in\R^d$.  Then
$P_*(\dd x)=Z(\theta_*)^{-1}\ee^{\langle\theta_*,x\rangle}\mu(\dd x)$ is a
probability measure with full affine support satisfying
\begin{equation}\label{eq:tilted-identity}
 \nu_{P_*}=C\ee^{\langle b,\theta_*\rangle}\,\delta_b*P_*^{*(d+1)}.
\end{equation}
\end{lemma}

\begin{proof}
Multiply both sides of \eqref{eq:measure-identity} by the positive function
$\ee^{\langle\theta_*,s\rangle}$.  On the left $s=x_0+\dots+x_d$, and on the
right $s=b+x_0+\dots+x_d$; this tilts every factor and produces the
constant $\ee^{\langle\theta_*,b\rangle}$.  Divide by $Z(\theta_*)^{d+1}$.
\end{proof}

Positivity of \eqref{eq:measure-identity} carries support information.
Let $S=\supp\mu$ and let $S^{\oplus m}$ be the set of sums of $m$ points of
$S$.  The weight in \eqref{eq:nu-def} is continuous and vanishes exactly on
affinely dependent tuples.  Hence $\supp\nu_\mu$ is the closure of the set
$\Sigma(S)$ of sums $x_0+\dots+x_d$ of affinely independent tuples in $S$,
and \eqref{eq:measure-identity} gives the \emph{support identity}
\begin{equation}\label{eq:support-identity}
 \overline{\Sigma(S)}=b+\overline{S^{\oplus(d+1)}}.
\end{equation}
Taking the infimum of a linear functional $q$ bounded below on $S$, with
$\alpha(q)=\inf_Sq$, gives
\begin{equation}\label{eq:basis-gap}
 \inf\Bigl\{\sum_{j=0}^d\bigl(q(x_j)-\alpha(q)\bigr):
 x_0,\dots,x_d\in S\text{ affinely independent}\Bigr\}
 =\langle q,b\rangle .
\end{equation}
In particular $\langle q,b\rangle\ge0$ for every such $q$, so $b$ is a
recession direction of $Q$, that is, $Q+tb\subset Q$ for $t\ge0$
\cite[Section~8]{Rockafellar1970}.

\section{Regularity and Constant Generalized Variance}
\label{sec:regularity}

This section proves \cref{thm:regularity-intro,thm:constant-intro} and
\cref{cor:jeffreys-intro,cor:regression-intro}.

\subsection{Regularity}

By \cref{lem:tilt}, a boundary point of $\Theta$ at which the Laplace
transform is finite would produce a probability measure that satisfies the
measure identity but whose Laplace transform is infinite at points
arbitrarily close to the origin.  We show that there is no such measure.
The argument uses the identity alone and not the exponential family, so we
state it for an arbitrary probability measure $P$ on $\R^d$ with full affine
support such that
\begin{equation}\label{eq:conv-identity}
 \nu_P*P^{*(n-d-1)}=C_0\,\delta_b*P^{*n}
\end{equation}
for some integer $n\ge d+1$, some $C_0>0$ and some $b\in\R^d$; here
$P^{*0}=\delta_0$, and $\nu_P$ is defined by \eqref{eq:nu-def} as a Borel
measure with values in $[0,\infty]$.  Since the right side has total mass
$C_0$, the identity forces $\nu_P(\R^d)=C_0$.  For $n=d+1$ it is
\eqref{eq:measure-identity}, the case needed for
\cref{thm:regularity-intro}; the additional convolution factors
arise from samples of size $n>d+1$ in
\cref{cor:regression-intro}.  The identity first gives rapid
one-dimensional tail decay.

\begin{lemma}[Quadratic Tail Recursion]\label{lem:tail}
Let $P$ be a probability measure with full affine support satisfying
\eqref{eq:conv-identity}.  For every nonzero linear
functional $q$ and $Y=\langle q,X\rangle$, $X\sim P$, there are
$A,R_0<\infty$ such that
\begin{equation}\label{eq:tail}
 P(Y>R)\le\frac A{R^2}\,P\Bigl(Y>\frac R{2n}\Bigr)\qquad(R\ge R_0).
\end{equation}
Consequently $P(Y>R)\le\exp\{-c(\log R)^2\}$ for some $c>0$ and all large
$R$, and $\mathbb E[Y_+^m]<\infty$ for every $m\ge1$.
\end{lemma}

\begin{proof}
Throughout the proof, $|\cdot|$ is the Euclidean norm, and $q$ is
identified with a vector of $\R^d$ through $\langle\cdot,\cdot\rangle$, so
that $|\langle q,x\rangle|\le|q|\,|x|$.  For $t\in\R$ let
$H_t=\{s\in\R^d:\langle q,s\rangle>t\}$.  We compare the two sides of
\eqref{eq:conv-identity} on $H_t$ for large $t$.  On the left, the weight
$\det(\widetilde x_0,\dots,\widetilde x_d)^2$ in \eqref{eq:nu-def} grows
like $|x_i|^2$ when one point $x_i$ is far out and the other points stay
near the vertices of a fixed simplex (Steps~1 and~2).  The right side is
$C_0$ times the law of $b$ plus a sum of $n$ independent random vectors
with law $P$ (Step~3).

\emph{Step 1: A Far Vertex Gives a Large Determinant.}  Since $P$ has full
affine support, $\supp P$ contains affinely independent points
$z_0,\dots,z_d$, and then
$\Delta=\det(\widetilde z_0,\dots,\widetilde z_d)\ne0$.  For $x\in\R^d$ let
$\beta_0(x),\dots,\beta_d(x)$ be the barycentric coordinates of $x$ with
respect to $z_0,\dots,z_d$, that is, the unique real numbers with
\[
 \sum_{i=0}^d\beta_i(x)=1,\qquad\sum_{i=0}^d\beta_i(x)\,z_i=x .
\]
These two conditions say that
$\sum_i\beta_i(x)\,\widetilde z_i=\widetilde x$, so Cramer's rule gives
\[
 \beta_i(x)=\frac{\det(\widetilde z_0,\dots,\widetilde z_{i-1},\widetilde x,
 \widetilde z_{i+1},\dots,\widetilde z_d)}{\Delta}\qquad(i=0,\dots,d).
\]
In particular each $\beta_i$ is affine; write
$\beta_i(x)=\beta_i(0)+\langle\ell_i,x\rangle$ with $\ell_i\in\R^d$.  If
$\langle\ell_i,v\rangle=0$ for every $i$, then
$v=\sum_i\beta_i(v)z_i=\sum_i\beta_i(0)z_i=0$.  Hence the continuous
function $v\mapsto\max_i|\langle\ell_i,v\rangle|$ is positive on the unit
sphere, and
\[
 \kappa=\min_{|v|=1}\,\max_{0\le i\le d}|\langle\ell_i,v\rangle|>0 .
\]
Geometrically, $|\beta_i(x)|$ is the distance from $x$ to the hyperplane
through the points $z_j$, $j\ne i$, divided by the distance from $z_i$ to
that hyperplane, and $\kappa>0$ says that a point far from the origin is
far from at least one of these $d+1$ hyperplanes.

We claim that there are $\rho>0$ and $R'<\infty$ with the following
property: for every $x\in\R^d$ with $|x|\ge R'$ there is an index $i$ such
that
\[
 |\det(\widetilde y_0,\dots,\widetilde y_d)|\ge c_1|x|,\qquad
 c_1=\tfrac14\,|\Delta|\,\kappa,
\]
for all $y_0,\dots,y_d\in\R^d$ with $y_i=x$ and $|y_j-z_j|<\rho$ for
$j\ne i$.  To prove the claim, fix $i$ and points $y_j$, $j\ne i$, and put
$y_i=x$.  The point $x$ enters $\det(\widetilde y_0,\dots,\widetilde y_d)$
only through the column $\widetilde y_i=\widetilde x$, so
$\det(\widetilde y_0,\dots,\widetilde y_d)/\Delta-\beta_i(x)$ is an affine
function of $x$.  Its constant term and its gradient are polynomials in the
$y_j$, $j\ne i$, and both vanish when $y_j=z_j$ for all $j\ne i$, by the
formula for $\beta_i$.  Choose $\rho>0$ such that, for every $i$, both have
norm at most $\kappa/2$ whenever $|y_j-z_j|<\rho$ for all $j\ne i$.  For
such $y_j$,
\[
 \Bigl|\frac{\det(\widetilde y_0,\dots,\widetilde y_d)}{\Delta}
 -\beta_i(x)\Bigr|\le\frac\kappa2\,(1+|x|).
\]
Now let $R'=2+4\max_j|\beta_j(0)|/\kappa$ and $|x|\ge R'$, and choose $i$
with $|\langle\ell_i,x\rangle|\ge\kappa|x|$, which is possible by the
definition of $\kappa$.  Then
$|\beta_i(x)|\ge\kappa|x|-|\beta_i(0)|$, so
\[
 \frac{|\det(\widetilde y_0,\dots,\widetilde y_d)|}{|\Delta|}
 \ge\kappa|x|-|\beta_i(0)|-\frac\kappa2\,(1+|x|)\ge\frac\kappa4\,|x| ,
\]
where the last step uses $|x|\ge R'$.  This proves the claim.

\emph{Step 2: A Lower Bound on Far Half-Spaces.}  Let
$U_j=\{y\in\R^d:|y-z_j|<\rho\}$, let $p_0=\min_jP(U_j)$, which is positive
because each $z_j$ lies in $\supp P$, and put
\[
 B=\sum_{j=0}^d\,\sup_{y\in U_j}|\langle q,y\rangle|,\qquad
 R_1=\max\Bigl\{|q|R',\ \max_{0\le j\le d}\,\sup_{y\in U_j}
 \langle q,y\rangle\Bigr\}.
\]
Let $R\ge R_1$, let $X_0,\dots,X_d$ be independent with law $P$, and for
$i=0,\dots,d$ let $E_i$ be the event that $X_i\in H_R$ and $X_j\in U_j$ for
all $j\ne i$.  By the choice of $B$, we have $X_0+\dots+X_d\in H_{R-B}$ on
$E_i$.  The events $E_i$ are pairwise disjoint: for $i'\ne i$, the event
$E_i$ requires $X_{i'}\in U_{i'}$, hence $\langle q,X_{i'}\rangle\le R$ by
the choice of $R_1$, while $E_{i'}$ requires $\langle q,X_{i'}\rangle>R$.
By the description of $\nu_P$ as an expectation in
\cref{sec:transform},
\begin{align*}
 (d+1)!\,\nu_P(H_{R-B})
 &=\mathbb E\bigl[\det(\widetilde X_0,\dots,\widetilde X_d)^2\,
 \mathbf 1_{\{X_0+\dots+X_d\in H_{R-B}\}}\bigr]\\
 &\ge\sum_{i=0}^d\mathbb E\bigl[\det(\widetilde X_0,\dots,
 \widetilde X_d)^2\,\mathbf 1_{E_i}\bigr].
\end{align*}
Fix $i$, and for $x\in\R^d$ put
\[
 g_i(x)=\mathbb E\bigl[\det(\widetilde X_0,\dots,\widetilde X_{i-1},
 \widetilde x,\widetilde X_{i+1},\dots,\widetilde X_d)^2\,
 \mathbf 1_{\{X_j\in U_j\text{ for all }j\ne i\}}\bigr].
\]
Since $X_i$ is independent of the other $X_j$, Fubini's theorem shows that
the $i$th term of the last sum is $\int_{H_R}g_i(x)\,P(\dd x)$.  Now let
$x\in H_R$.  Then $|x|\ge\langle q,x\rangle/|q|>R/|q|\ge R'$, so Step~1
provides an index $i$ for which the determinant in the definition of
$g_i(x)$ has absolute value at least $c_1|x|$ on the event
$\{X_j\in U_j\text{ for all }j\ne i\}$.  This event has probability at
least $p_0^d$, so $\sum_ig_i(x)\ge c_1^2p_0^d|x|^2$.  Since
$|x|^2\ge\langle q,x\rangle^2/|q|^2$, integrating over $H_R$ gives
\[
 \nu_P(H_{R-B})\ge c_2\,\mathbb E[Y^2\mathbf 1_{\{Y>R\}}]\quad(R\ge R_1),
 \qquad c_2=\frac{c_1^2p_0^d}{(d+1)!\,|q|^2}.
\]

\emph{Step 3: An Upper Bound from the Identity.}  Choose $B_1>0$ such that
$p_1=P^{*(n-d-1)}(H_{-B_1})$ is positive, which is possible because
$H_{-B_1}$ increases to $\R^d$ as $B_1\to\infty$; if $n=d+1$, then
$P^{*0}=\delta_0$ and $p_1=1$.  If $s\in H_{R-B}$ and $y\in H_{-B_1}$, then
$s+y\in H_{R-B-B_1}$.  Hence, using \eqref{eq:conv-identity} for the
equality,
\[
 p_1\,\nu_P(H_{R-B})\le\bigl(\nu_P*P^{*(n-d-1)}\bigr)(H_{R-B-B_1})
 =C_0\,\bigl(\delta_b*P^{*n}\bigr)(H_{R-B-B_1}).
\]
The right side is $C_0\,P(Y_1+\dots+Y_n>r)$, where
$r=R-B-B_1-\langle q,b\rangle$ and $Y_1,\dots,Y_n$ are independent copies
of $Y$.  If a sum of $n$ numbers exceeds $r$, one summand exceeds $r/n$.
When $R\ge2(B+B_1+|\langle q,b\rangle|)$, we have $r\ge R/2$, and the union
bound gives
\[
 P(Y_1+\dots+Y_n>r)\le n\,P(Y>r/n)\le n\,P\Bigl(Y>\frac R{2n}\Bigr).
\]

\emph{Step 4: The Recursion.}  Let $R\ge R_0$, where
$R_0=\max\{R_1,\,2(B+B_1+|\langle q,b\rangle|)\}$.  Since $Y^2>R^2$ on
$\{Y>R\}$, Steps~2 and~3 give
\[
 c_2R^2P(Y>R)\le c_2\,\mathbb E[Y^2\mathbf 1_{\{Y>R\}}]
 \le\nu_P(H_{R-B})
 \le\frac{nC_0}{p_1}\,P\Bigl(Y>\frac R{2n}\Bigr).
\]
This is \eqref{eq:tail} with $A=nC_0/(c_2p_1)$.  Every
term lies in $[0,\infty]$, so no moment of $P$ has been assumed.

\emph{Step 5: Iteration and Moments.}  Fix $R\ge R_0$ and let
$k=1+\lfloor\log(R/R_0)/\log(2n)\rfloor$.  Then $R/(2n)^j\ge R_0$ for
$0\le j\le k-1$.  Applying \eqref{eq:tail} at
$R,R/(2n),\dots,R/(2n)^{k-1}$ and bounding the last probability by $1$,
we get
\[
 P(Y>R)\le\prod_{j=0}^{k-1}\frac{A\,(2n)^{2j}}{R^2}
 =A^k(2n)^{k(k-1)}R^{-2k}.
\]
Since $k=\log R/\log(2n)+O(1)$, the logarithm of the right side is
\[
 k\log A+k(k-1)\log(2n)-2k\log R
 =-\frac{(\log R)^2}{\log(2n)}+O(\log R).
\]
Hence, for every $c<1/\log(2n)$ and all large $R$,
\[
 P(Y>R)\le\exp\{-c(\log R)^2\}=R^{-c\log R}.
\]
The right side decays faster than every power of $R$, so
\[
 \mathbb E[Y_+^m]=m\int_0^\infty R^{m-1}P(Y>R)\,\dd R<\infty
 \qquad(m\ge1).\qedhere
\]
\end{proof}

The moments supplied by \cref{lem:tail} make the characteristic
function of $P$ smooth.  Fourier transformation then turns
\eqref{eq:conv-identity} into an analytic elliptic equation for its
logarithm, and analyticity yields exponential moments.

\begin{proposition}[Exponential Moments]\label{prop:exp-moments}
Let $P$ be a probability measure on $\R^d$ with full affine support
satisfying \eqref{eq:conv-identity}.  Then
$\int\ee^{\varepsilon|x|}\,P(\dd x)<\infty$ for some $\varepsilon>0$.
Consequently $P$ generates a natural exponential family whose natural
parameter domain contains a neighbourhood of the origin, and this family
satisfies \eqref{eq:E} with $C=C_0$.
\end{proposition}

\begin{proof}
\Cref{lem:tail}, applied to each coordinate functional and its
negative, gives $\mathbb E_P|X|^m<\infty$ for all $m$.  Hence the
characteristic function
$\phi(t)=\mathbb E_P\ee^{i\langle t,X\rangle}$ is $C^\infty$.  On a ball
where $\phi\neq0$, let $L(t)=\log\phi(t)$ with $L(0)=0$.  The Cauchy--Binet
computation of \cref{lem:KS}, applied to the complex weight
$\ee^{i\langle t,x\rangle}$, shows that the Fourier transform of $\nu_P$ is
$(-1)^d\phi^{d+1}\det D^2L$ on this ball.  Fourier transformation of
\eqref{eq:conv-identity} and division by $\phi^n$ therefore give
\begin{equation}\label{eq:fourier-MA}
 \det D^2L(t)=(-1)^dC_0\,\ee^{i\langle b,t\rangle}.
\end{equation}
All differentiations are justified by the moment bounds.  Write $L=U+iV$
and
\begin{equation*}
 a(\xi)=\cof(D^2L)^{ij}\xi_i\xi_j.
\end{equation*}
Since
$D^2L(0)=-\Cov_P(X)$ is negative definite,
$\operatorname{Re}\,((-1)^{d-1}a(\xi))>0$ for $\xi\ne0$, uniformly on a
small ball.  On the pair $(U,V)$, the linearization of \eqref{eq:fourier-MA}
has symbol matrix
$\bigl(\begin{smallmatrix}\operatorname{Re}a&-\operatorname{Im}a\\
\operatorname{Im}a&\operatorname{Re}a\end{smallmatrix}\bigr)$, whose
quadratic form is $\operatorname{Re}a(\xi)(\eta_1^2+\eta_2^2)$ up to the
harmless sign $(-1)^{d-1}$.  Thus \eqref{eq:fourier-MA} is a strongly
elliptic analytic nonlinear system near $0$.  Morrey proved the following
\cite[p.~198]{Morrey1958}: let $u=(u^k)$ solve an analytic nonlinear
system whose $j$th equation involves derivatives of $u^k$ of order at most
$s_j+t_k$ and which is elliptic along $u$ in the sense of Douglis and
Nirenberg; if each $u^k$ has H\"older continuous derivatives of order
$t_k$, then $u$ is analytic at each interior point of its domain.  Here
$u=(U,V)$ is smooth, $s_j=0$ and $t_k=2$, so $\phi$ is real analytic near
$0$.  For $d=1$ the equation reads $L''=-C_0\ee^{ibt}$ and is integrated
directly.

Real analyticity gives $M,r>0$ with $\mathbb E_P\langle q,X\rangle^{2k}\le
M(2k)!\,r^{-2k}$ for all unit $q$ and all $k\ge0$.  By monotone convergence
in the series for $\cosh$,
$\mathbb E_P\exp(r'|\langle q,X\rangle|)<\infty$ for $0<r'<r$, and
H\"older's inequality gives the first assertion.

Consequently the Laplace transform $Z$ of $P$ is finite near the origin,
and $P$ generates an NEF with $0\in\Theta$.  For $\theta\in\Theta$,
\cref{lem:KS} shows that the Laplace transform of the left side of
\eqref{eq:conv-identity} is $Z(\theta)^n\det\nabla^2K(\theta)$, while that
of the right side is $C_0\ee^{\langle b,\theta\rangle}Z(\theta)^n$.  This is
\eqref{eq:E} with $C=C_0$.
\end{proof}

Conversely, by \cref{lem:KS,lem:tilt}, every member of an
NEF that satisfies \eqref{eq:E} satisfies \eqref{eq:conv-identity} with
$n=d+1$.  Thus the probability measures with full affine support that
satisfy \eqref{eq:conv-identity} are exactly the members of the natural
exponential families that satisfy \eqref{eq:E}; the finiteness of the
Laplace transform on an open set is a consequence of the identity and not
an additional assumption.  \Cref{thm:regularity-intro} is now
immediate.

\begin{proof}[Proof of \cref{thm:regularity-intro}]
Suppose that $\theta_*\in\partial\Theta$ and $Z(\theta_*)<\infty$.  By
\cref{lem:tilt}, the tilt $P_*$ of $\mu$ at $\theta_*$ is a
probability measure with full affine support satisfying
\eqref{eq:conv-identity} with $n=d+1$.  By
\cref{prop:exp-moments} it has exponential moments in a full
neighbourhood of the origin, which puts $\theta_*$ in the interior of $D$,
that is, in $\Theta$.  This is a contradiction.

Hence $D=\Theta$.  If $\theta_j\to\theta_*\in\partial\Theta$, Fatou's
lemma gives $\liminf Z(\theta_j)\ge Z(\theta_*)=\infty$, so
$K(\theta_j)\to\infty$.  An NEF with open natural parameter domain is
steep, and a steep NEF has mean domain $\inter Q$
\cite[Chapters~8--9]{BarndorffNielsen1978}.
\end{proof}

\subsection{The Constant Case}

With regularity in hand, the constant case reduces to the following form of
the J\"orgens--Calabi--Pogorelov theorem, due to Cheng and Yau: a smooth
convex solution of $\det D^2u=1$ on a convex domain whose graph is closed in
$\R^{d+1}$ is a quadratic polynomial
\cite{Jorgens1954,Calabi1958,Pogorelov1972,ChengYau1986}; see also
\cite[Theorem~1]{Loftin2008}.

\begin{proof}[Proof of \cref{thm:constant-intro}]
Here $b=0$.  By \cref{thm:regularity-intro}, $K(\theta)\to+\infty$
at every finite boundary point of the convex open set $\Theta$, so the
graph of $u=C^{-1/d}K$ over $\Theta$ is a closed subset of $\R^{d+1}$.
The function $u$ is smooth, strictly convex and satisfies $\det D^2u=1$.
By the theorem just quoted, $u$ is a quadratic polynomial.  A
quadratic polynomial is finite everywhere, so $\Theta=\R^d$, and
$K(\theta)=\frac12\theta^{\mathsf T}\Sigma\theta+\langle m,\theta\rangle+k$
with $\Sigma>0$.  Uniqueness of Laplace transforms gives
$\mu=\ee^k\,N(m,\Sigma)$.
\end{proof}

\begin{proof}[Proof of \cref{cor:jeffreys-intro}]
The Jeffreys density is proportional to $(\det\nabla^2K)^{1/2}$; if it is
constant, so is the generalized variance.
\end{proof}

For \cref{cor:regression-intro} we express the regression
condition as the measure identity \eqref{eq:conv-identity} with $b=0$ and
apply \cref{prop:exp-moments}.

\begin{proof}[Proof of \cref{cor:regression-intro}]
Let $P$ be the common law of the $X_j$ and put $S=X_1+\dots+X_n$.  For a
subset $J$ of $\{1,\dots,n\}$ with $d+1$ elements, let $\Delta_J$ be the
determinant of the vectors $\widetilde X_j$, $j\in J$.  The Schur
complement used in the proof of \cref{lem:KS} and the Cauchy--Binet
formula give
\begin{equation}\label{eq:scatter-det}
 n\det W
 =\det\Bigl(\sum_{j=1}^n\widetilde X_j\widetilde X_j^{\mathsf T}\Bigr)
 =\sum_J\Delta_J^2 .
\end{equation}
Let $f\ge0$ be a bounded Borel function on $\R^d$.  By exchangeability the
$\binom n{d+1}$ expectations $\mathbb E[\Delta_J^2f(S)]$ are equal, and by
\eqref{eq:nu-def} each of them equals
$(d+1)!\int f\,\dd(\nu_P*P^{*(n-d-1)})$.  Hence
\begin{equation}\label{eq:regression-measure}
 \mathbb E[\det W\,f(S)]
 =\frac{(n-1)!}{(n-d-1)!}\int f\,\dd\bigl(\nu_P*P^{*(n-d-1)}\bigr).
\end{equation}
Since $\bar X=S/n$, it follows that
$\mathbb E[\det W\mid\bar X]=c$ almost surely, for a constant
$c\in(0,\infty)$, if and only if
\begin{equation}\label{eq:regression-identity}
 \nu_P*P^{*(n-d-1)}=C_0\,P^{*n},
 \qquad C_0=\frac{(n-d-1)!}{(n-1)!}\,c .
\end{equation}

Suppose that \eqref{eq:regression-identity} holds.  Then $\det W$ is not
almost surely zero, so $P$ has full affine support, and
\eqref{eq:conv-identity} holds with $b=0$.  By
\cref{prop:exp-moments}, $P$ generates an NEF with constant
generalized variance, and \cref{thm:constant-intro} shows that $P$
is a nondegenerate normal distribution.  Conversely, let $P=N(m,\Sigma)$
with $\Sigma>0$.  Its cumulant function is quadratic with Hessian $\Sigma$,
so \cref{lem:KS} gives $\nu_P=\det\Sigma\,P^{*(d+1)}$.  Convolution
with $P^{*(n-d-1)}$ gives \eqref{eq:regression-identity} with
$C_0=\det\Sigma$, that is, with $c=\det\Sigma\,(n-1)!/(n-d-1)!$.
\end{proof}

\begin{remark}[Regularity Cannot Be Dropped]\label{rem:not-cumulant}
The function $u(x,y)=x^2/(2y)+y^3/6$ is strictly convex on $\{y>0\}$ and
satisfies $\det D^2u=1$, but $u(0,y)\to0$ as $y\downarrow0$, so its graph
is not closed.  By \cref{thm:regularity-intro}, no positive multiple
of $u$ plus an affine function is the cumulant function of a positive
measure.  Positivity of $\mu$ enters \cref{thm:constant-intro} only
through regularity.
\end{remark}

\section{Boundary Structure}
\label{sec:boundary}

This section proves \cref{thm:boundary-intro}.  Assume that $b\ne0$.
The proof has three steps.  The Legendre dual of \eqref{eq:E} yields a
sequence of means that approaches the boundary of the mean domain with
bounded relative entropy.  Along such a sequence the laws $P_\theta$
concentrate near an exposed face of $Q$, which $\mu$ must therefore charge.
The support identity \eqref{eq:support-identity} then forces the lattice
structure of $\mu$ in the directions normal to that face.

\subsection{The Dual Potential}

Normalize $0\in\Theta$ and $Z(0)=1$, so that $\mu=P_0$ is a probability
measure, and let $u=K^*$ be the Legendre transform of $K$.  It is a smooth
strictly convex function on the mean domain $M$, which we call the
\emph{dual potential}.  Then
$\nabla u(m)=\theta$ and $\nabla^2u(m)=\nabla^2K(\theta)^{-1}$ for
$m=\nabla K(\theta)$, so \eqref{eq:E} becomes
\begin{equation}\label{eq:dual}
 \det D^2u=\exp\{-\langle b,\nabla u\rangle-\log C\}\qquad\text{on }M,
\end{equation}
the affine K\"ahler--Ricci-flat equation of~\cite{LiXu2009}.  Moreover
$u(m)=\KL(P_\theta\Vert\mu)$, where
$\KL(P\Vert R)=\int\log(\dd P/\dd R)\,\dd P$ denotes relative entropy, and
more generally, for $\theta_0\in\Theta$ with mean $m_0$,
\begin{equation}\label{eq:bregman}
 u(m)-u(m_0)-\langle\theta_0,m-m_0\rangle=\KL(P_\theta\Vert P_{\theta_0}).
\end{equation}

Sublevel sets of this relative entropy are bounded in the mean parameter.

\begin{lemma}[Bounded Entropy Bounds the Mean]\label{lem:entropy-mean}
For every $\theta_0\in\Theta$ and $A<\infty$, the set
$\{m\in M:\KL(P_\theta\Vert P_{\theta_0})\le A\}$ is bounded.
\end{lemma}

\begin{proof}
Choose $\varepsilon>0$ with $\theta_0+\overline{B(0,\varepsilon)}\subset
\Theta$.  For $m\ne m_0$ let $v=(m-m_0)/|m-m_0|$.  The Donsker--Varadhan
variational inequality
$\mathbb E_Pf\le\KL(P\Vert R)+\log\mathbb E_R\ee^f$ gives
\begin{align*}
 \varepsilon\langle v,m-m_0\rangle
 &\le\KL(P_\theta\Vert P_{\theta_0})
   +\log \mathbb E_{P_{\theta_0}}\ee^{\varepsilon\langle v,X-m_0\rangle}\\
 &\le A+\max_{|w|=1}\log
   \mathbb E_{P_{\theta_0}}\ee^{\varepsilon\langle w,X-m_0\rangle}.
 \qedhere
\end{align*}
\end{proof}

The next proposition is the Liouville-type input suggested in
\cite[arXiv version, Remark~2]{LiXu2009}.  We adapt Wang's proof for entire
solutions \cite[Lemma~2.3, Proposition~2.4 and the proof of
Proposition~2.2]{Wang2026}, whose auxiliary function in Step~1 goes back to
Wu and Xu~\cite{WuXu2021}; on a convex domain it uses only the compactness
of sublevel sets.

\begin{proposition}[Drifted Monge--Amp\`ere Nonexistence]
\label{prop:dual-nonexistence}
Let $U\subset\R^d$ be a convex open set, $a\in\R^d\setminus\{0\}$,
$c\in\R$, and let $u\in C^4(U)$ satisfy $D^2u>0$ and
\begin{equation}\label{eq:drift-MA}
 \det D^2u=\exp\{-\langle a,\nabla u\rangle-c\}\qquad\text{on }U.
\end{equation}
Suppose that for every $y\in U$ and $t\in\R$ the set
$\{x\in U:u_y(x)\le t\}$ is compact, where
$u_y(x)=u(x)-u(y)-\langle\nabla u(y),x-y\rangle$.  Then no such $u$ exists.
\end{proposition}

\begin{proof}
First let $d=1$ and $w=u'$.  Then $w'=\ee^{-aw-c}$, so
$(\ee^{aw})'=a\ee^{-c}$ and $\ee^{aw(x)}=a\ee^{-c}(x-x_0)$ for some $x_0$.
Thus $U$ lies in the half-line $\{a(x-x_0)>0\}$ and has a finite endpoint.
Near it, $w$ has at most a logarithmic singularity and is integrable, so
$u_y$ stays bounded, and a sublevel set of $u_y$ accumulates at a point
outside $U$.  This contradicts compactness.

Since $U$ is convex and $D^2u>0$, each $u_y$ is nonnegative.  Let $d\ge2$.
\emph{Step 1: A Gradient Bound.}  We show that
\begin{equation*}
 |\langle a,\nabla u(x)-\nabla u(y)\rangle|\le 3d+3\quad\text{on } U.
\end{equation*}
Fix $y$, put $v=u_y+3$, and let $c_y=c+\langle
a,\nabla u(y)\rangle$.  Then $v(y)=3$, $v\ge3$, and $v$ satisfies
\eqref{eq:drift-MA} with $c_y$ in place of $c$.  For $T>3$ the set
$U_T=\{v<T\}$ has compact closure in $U$.  Put $\eta=T-v$, $\kappa=3T$ and
\[
 G=\exp\{-\kappa/\eta-\langle a,\nabla v\rangle-c_y\}
 =\ee^{-\kappa/\eta}\det D^2v\qquad\text{on }U_T .
\]
$G$ extends continuously by zero to $\partial U_T$ and attains its maximum
at some $x_T\in U_T$.  Write $(v^{ij})=(D^2v)^{-1}$.  At $x_T$,
$0=\partial_i\log G=-\kappa v_i/\eta^2-a_\ell v_{\ell i}$.  Contracting with
$a_i$ and with $v^{ij}v_j$ gives
\begin{equation}\label{eq:wang-first}
 a^{\mathsf T}D^2v\,a=\frac\kappa{\eta^2}X,
 \qquad
 X:=-\langle a,\nabla v(x_T)\rangle=\frac\kappa{\eta^2}v^{ij}v_iv_j\ge0 .
\end{equation}
Differentiating the logarithm of \eqref{eq:drift-MA} gives
$v^{ij}v_{ij\ell}=-a_kv_{k\ell}$.  Hence, at $x_T$,
\[
 0\ge v^{ij}\partial_{ij}\log G
 =a^{\mathsf T}D^2v\,a-\frac{\kappa d}{\eta^2}
 -\frac{2\kappa}{\eta^3}v^{ij}v_iv_j
 =\Bigl(\frac\kappa{\eta^2}-\frac2\eta\Bigr)X-\frac{\kappa d}{\eta^2}.
\]
Since $\kappa-2\eta=T+2v>T$ on $U_T$, this gives $X\le\kappa d/(\kappa-2\eta)
<3d$.  For $x$ with $v(x)<T$, the inequality $G(x)\le G(x_T)$ gives
$-\langle a,\nabla v(x)\rangle\le X+\kappa/(T-v(x))<3d+3T/(T-v(x))$.
Letting $T\to\infty$ yields $\langle a,\nabla u(y)-\nabla u(x)\rangle\le3d+3$.
Exchanging $x$ and $y$ proves Step~1.

\emph{Step 2: The Maximum-Principle Contradiction.}  Fix $y\in U$ and
replace $u$ by $u_y$.  This preserves
\eqref{eq:drift-MA} with a new constant, as well as the hypothesis.  Now
$u\ge0$, $u(y)=0$ and $\nabla u(y)=0$, and by Step~1,
$w=\langle a,\nabla u\rangle$ satisfies $|w|\le W:=3d+3$.  Let
$g_{ij}=u_{ij}$, $(g^{ij})=(g_{ij})^{-1}$, and put
\[
 P=a^{\mathsf T}D^2u\,a>0,\qquad T_{ij}=a_ku_{ijk},\qquad
 \mathcal Lf=g^{ij}f_{ij}+a_kf_k .
\]
Differentiating $\log\det D^2u=-\langle a,\nabla u\rangle-c$ once and twice
gives
\begin{equation*}
 g^{ij}u_{ijk}=-a_\ell u_{\ell k},\qquad
 g^{ij}u_{ijrs}-g^{ip}g^{jq}u_{ijr}u_{pqs}=-a_\ell u_{\ell rs}.
\end{equation*}
Contracting with $a_k$, respectively $a_ra_s$, and using
$P_i=T_{ij}a_j$, we get
\begin{equation}\label{eq:wang-identities}
 g^{ij}T_{ij}=-P,\qquad \mathcal LP=|T|_g^2 .
\end{equation}
At a given point, choose a $g$-orthonormal basis whose first vector is
parallel to $a$, so that $a=\sqrt P\,e_1$, and write $\tau=T_{11}$.  Then
$\operatorname{tr}T=-P$, $|\nabla P|_g^2=P\sum_iT_{1i}^2$ and
$\langle a,\nabla P\rangle=P\tau$.  Let $\vartheta=1/(4(d-1))$.  Then
\[
 \mathcal L(-P^{-\vartheta})
 =\vartheta P^{-\vartheta-1}\Bigl(|T|^2-(\vartheta+1)\sum_iT_{1i}^2\Bigr).
\]
Since $T$ is symmetric, it has the block form
$T=\bigl(\begin{smallmatrix}\tau&v^{\mathsf T}\\ v&B\end{smallmatrix}\bigr)$.  Then
\begin{gather*}
 |T|^2-(\vartheta+1)\sum_iT_{1i}^2=-\vartheta\tau^2+(1-\vartheta)|v|^2+|B|^2,\\
 |B|^2\ge(\operatorname{tr}B)^2/(d-1)=(P+\tau)^2/(d-1).
\end{gather*}
Hence
\begin{equation}\label{eq:wang-subsolution}
 \mathcal L(-P^{-\vartheta})\ge\vartheta P^{-\vartheta-1}
 \Bigl(\frac{(P+\tau)^2}{d-1}-\vartheta\tau^2\Bigr).
\end{equation}
For $\varepsilon>0$ put $F_\varepsilon=-P^{-\vartheta}-\varepsilon u\le
-\varepsilon u$.  Its superlevel sets are compact, so it attains its
maximum at some $x_\varepsilon\in U$.  There $\nabla F_\varepsilon=0$ and
$D^2F_\varepsilon\le0$, so $\mathcal LF_\varepsilon(x_\varepsilon)\le0$ and
$\vartheta P^{-\vartheta-1}\nabla P=\varepsilon\nabla u$.  Pairing the last
identity with $a$ gives $\vartheta P^{-\vartheta}\tau=\varepsilon w$, hence
$|\tau|\le\varepsilon WP^\vartheta/\vartheta$.

If $|\tau|\ge P/2$, then $P^{1-\vartheta}\le2\varepsilon W/\vartheta$.  If
$|\tau|<P/2$, the bracket in \eqref{eq:wang-subsolution} is at least
\begin{equation*}
 P^2/(4(d-1))-\vartheta P^2/4=3P^2/(16(d-1)).
\end{equation*}
Since $\mathcal Lu=d+w\le d+W$, this gives
\begin{equation*}
 0\ge\mathcal LF_\varepsilon\ge3\vartheta P^{1-\vartheta}/(16(d-1))-\varepsilon(d+W).
\end{equation*}
In either case
$P(x_\varepsilon)\to0$ as $\varepsilon\downarrow0$.  But $F_\varepsilon(x_\varepsilon)\ge
F_\varepsilon(y)=-P(y)^{-\vartheta}$ and $F_\varepsilon\le-P^{-\vartheta}$, so
$P(x_\varepsilon)\ge P(y)>0$, a contradiction.
\end{proof}

Applied to the dual potential, \cref{prop:dual-nonexistence}
shows that some sublevel set of the entropy reaches the boundary of the mean
domain.

\begin{theorem}[Finite-Entropy Boundary Access]\label{thm:access}
Let $\mu$ satisfy \eqref{eq:E} with $b\ne0$.  Then there are $A<\infty$ and
means $m_n\in M$ with $u(m_n)\le A$ and $m_n\to m_*\in\partial M$.
\end{theorem}

\begin{proof}
Suppose not, and let $\theta_0\in\Theta$ be arbitrary.  By
\eqref{eq:bregman}, the
function $u_{m_0}$ of \cref{prop:dual-nonexistence} is
$\KL(P_\theta\Vert P_{\theta_0})$, so its sublevel sets are bounded by
\cref{lem:entropy-mean}.  The functions $u_{m_0}$ and $u$ differ by an
affine function, which is bounded on bounded sets.  Hence no sublevel set of
$u_{m_0}$ accumulates at a finite point of $\partial M$, and all sublevel
sets are compact subsets of $M$.  Since $M=\inter Q$ is open and convex by
\cref{thm:regularity-intro} and $u$ is smooth with $D^2u>0$,
\cref{prop:dual-nonexistence},
applied to \eqref{eq:dual} with $a=b$, gives a contradiction.
\end{proof}

Along such a sequence the natural parameters diverge, and their limiting
direction exposes a face of $Q$ that $\mu$ must charge.

\begin{theorem}[Positive-Mass Exposed Face]\label{thm:massive-face}
Let $(m_n)\subset M$ converge to $m_*\in\partial M$, and suppose that
$u(m_n)\le A<\infty$.  Put $\theta_n=\nabla u(m_n)$.  Then
$|\theta_n|\to\infty$.  If $\theta_n/|\theta_n|\to v$ along a subsequence,
the exposed face
$F=\{x\in Q:\langle v,x\rangle=h_Q(v)\}$, where
$h_Q(v)=\sup_Q\langle v,\cdot\rangle$ is the support function of $Q$, is
proper, contains $m_*$, and satisfies $\mu(F)>0$.
\end{theorem}

\begin{proof}
If a subsequence of $(\theta_n)$ were bounded, it would have a limit point
$\theta_*$.  Since
\begin{equation*}
 K(\theta_n)=\langle\theta_n,m_n\rangle-u(m_n)
\end{equation*}
is bounded above, Fatou's lemma gives $Z(\theta_*)<\infty$.  Hence
$\theta_*\in\Theta$ by \cref{thm:regularity-intro}, and
$m_*=\nabla K(\theta_*)\in M$, a contradiction.

So $|\theta_n|\to\infty$.  Let $r_n=|\theta_n|$ and $v_n=\theta_n/r_n\to v$.
Dividing $K(\theta_n)=\langle\theta_n,m_n\rangle-u(m_n)$ by $r_n$ gives
$K(\theta_n)/r_n\to\langle v,m_*\rangle$.  For each support point $x$ and
$\delta>0$, $Z(\theta_n)\ge\mu(B(x,\delta))\ee^{r_n(\langle
v_n,x\rangle-\delta)}$, whence $h_Q(v)\le\langle v,m_*\rangle$.  The reverse
inequality holds because $m_*\in Q$, so $m_*\in F$; and $F$ is proper
because $v\ne0$ and $\supp\mu$ spans $\R^d$ affinely.  Put
$Y=h_Q(v)-\langle v,X\rangle\ge0$.  Then $\mathbb E_{P_{\theta_n}}Y\to0$, so Markov's
inequality gives $P_{\theta_n}(E_\varepsilon)\to1$ for
$E_\varepsilon=\{Y\le\varepsilon\}$.  Data processing for relative entropy
through the partition $\{E_\varepsilon,E_\varepsilon^c\}$ yields
$\liminf_n\KL(P_{\theta_n}\Vert\mu)\ge-\log\mu(E_\varepsilon)$.  If
$\mu(F)=0$, then $\mu(E_\varepsilon)\downarrow0$ as $\varepsilon\downarrow0$,
contradicting $u(m_n)\le A$.
\end{proof}

The second half of this proof is a special case of
\cite[Theorem~1(ii)]{CsiszarMatus2005} combined with uniform integrability
of bounded-entropy densities.  \Cref{thm:access} supplies a sequence
to which \cref{thm:massive-face} applies, and together they give the
following.

\begin{corollary}[Existence of a Positive-Mass Face]\label{cor:face-null}
Let $\mu$ satisfy \eqref{eq:E}.
\begin{enumerate}
 \item[\textup{(i)}] If $\mu$ is not Gaussian, some proper exposed face $F$
 of $Q$ has $\mu(F)>0$.  In particular $\mu$ is Gaussian if it is absolutely
 continuous or if $Q=\R^d$.
 \item[\textup{(ii)}] $\supp\mu$ is unbounded.
\end{enumerate}
\end{corollary}

\begin{proof}
For (i): if $b=0$, use \cref{thm:constant-intro}; if $b\ne0$, use
\cref{thm:access,thm:massive-face}.  Proper exposed faces
of $Q$ lie in hyperplanes, and $Q=\R^d$ has no proper exposed face.  For
(ii): bounded support gives $D=\R^d$ and
$\det\nabla^2K\le(\operatorname{diam}\supp\mu)^{2d}$, so $b=0$ (let
$\theta=tb$, $t\to\infty$).  Then $\mu$ is Gaussian by
\cref{thm:constant-intro}, which contradicts bounded support.
\end{proof}

\subsection{Normal Structure at a Positive-Mass Face}

For the rest of this section, $F$ is a proper exposed face with $\mu(F)>0$,
and we describe $\mu$ in the directions normal to $F$.  Translate a point of
$\supp(\mu|_F)$ to the origin and put
\[
 V=\Span(\supp\mu|_F),\quad r=\dim V,\quad N=\R^d/V,\quad k=d-r,
\]
let $\pi:\R^d\to N$ be the quotient map, and put $\bar\mu=\pi_\#\mu$ and
$\beta=\pi(b)$.
Points of $\supp\mu$ in $V$ lie in $F$, so $\bar\mu$ has an atom at $0$ of
mass $\mu(F)$.  Choose an affine function $h=\langle a,\cdot\rangle-\alpha$
with $h\ge0$ on $Q$ and $F=Q\cap\{h=0\}$; we call $h(x)$ the \emph{height}
of $x$.  The function $h$ vanishes on the affine hull of $F$ and therefore
induces a functional $\bar h$ on $N$.

The first step produces a basis of $N$ consisting of atoms of $\bar\mu$.

\begin{lemma}[Minimal Normal Basis]\label{lem:normal-basis}
There are nonzero atoms $g_1,\dots,g_k$ of $\bar\mu$ forming a basis of $N$
with $\beta=g_1+\dots+g_k$.  Moreover $\langle a,b\rangle>0$.
\end{lemma}

\begin{proof}
Use the normal tilting in \cite[Proposition~3.4]{Masmoudi2007}.  Since
$\ee^{-th}\le1$ on $\supp\mu$ and $\int|x|^2\mu(\dd x)<\infty$, dominated
convergence gives $P_{-ta}\to\mu(F)^{-1}\mu|_F$ in total variation and in
the first two moments as $t\to\infty$.  The covariance limit is singular,
so $\det\nabla^2K(-ta)\to0$, and \eqref{eq:E} forces
$\langle a,b\rangle>0$.

Project \eqref{eq:measure-identity} to $N$.  The right side has an atom at
$\beta$ of mass at least $C\mu(F)^{d+1}$, so the projected left side has an
atom there as well.  By Fubini's theorem, a set of positive
$\det^2\cdot\mu^{\otimes(d+1)}$-measure of tuples with
$\sum_j\pi(x_j)=\beta$ consists of tuples whose projections are atoms of
$\bar\mu$.  Therefore there are affinely independent
$x_0,\dots,x_d\in\supp\mu$ with $\pi(x_j)$ atoms of $\bar\mu$ and
$\sum_j\pi(x_j)=\beta$.  Their total height is
$\sum_jh(x_j)=\bar h(\beta)=\langle a,b\rangle$, which by
\eqref{eq:basis-gap} is the minimum over all affinely independent tuples in
$\supp\mu$.

If $h(x_j)=0$, then $y=\pi(x_j)$ satisfies $\bar h(y)=0$, so
$\pi^{-1}(y)\cap Q\subset F$ and
$\bar\mu\{y\}=\mu|_F(\pi^{-1}(y))$.  This vanishes unless $y=0$, because
$\mu|_F$ lives on $V=\pi^{-1}(0)$.  Hence the points of height zero in the
tuple lie in $V$, and there are at most $r+1$ of them.  Suppose they did not
span $V$ affinely.  Then some $z\in\supp(\mu|_F)$ would lie outside their
affine span, and in the expansion of $\widetilde z$ in the basis
$\{\widetilde x_j\}$ of $\R^{d+1}$ some point of positive height would have a nonzero
coefficient.  Exchanging that point for $z$ would give an affinely
independent tuple of smaller total height, a contradiction.  Hence the tuple
consists of $r+1$ points spanning $V$ and $k$ points of positive height.
The projections $g_1,\dots,g_k$ of the latter form a basis of $N$, and they
sum to $\beta$ because the points of $V$ project to $0$.
\end{proof}

The support of $\bar\mu$ lies in the monoid $\Gamma$ generated by these
atoms.  Since they are linearly independent, $\Gamma$ is a free commutative
monoid.

\begin{theorem}[Free Commutative Monoid]\label{thm:free-monoid}
With $g_1,\dots,g_k$ as in \cref{lem:normal-basis},
\begin{equation}\label{eq:free-monoid}
 \supp\bar\mu\subset\Gamma:=\Nzero g_1+\dots+\Nzero g_k .
\end{equation}
\end{theorem}

\begin{proof}
Let $\bar S=\supp\bar\mu$.  If $y\in\bar S$ has nonzero $i$th coordinate in
the basis $(g_j)$, then $y$ can replace $g_i$ in a normal basis.  Take
$r+1$ points spanning $V$ from $\supp(\mu|_F)$, points of $\supp\mu$ over
the $g_j$ with $j\ne i$, and a point of $\supp\mu$ whose projection is close
to $y$.  Together they form an affinely independent tuple.  The projected support identity
\eqref{eq:support-identity} then gives
\begin{equation}\label{eq:descent}
 y-g_i\in\overline{\bar S^{\oplus(d+1)}}.
\end{equation}
Let $\delta=\min_i\bar h(g_i)>0$.  If $0\ne y\in\bar S$, some coordinate of $y$
is nonzero, and \eqref{eq:descent} gives $\bar h(y)\ge\bar h(g_i)\ge\delta$,
because $\bar h\ge0$ on $\bar S$.  So $\bar S\cap\{\bar h<\delta\}=\{0\}$.
Suppose that every point of $\bar S$ of height less than $n\delta$ lies in
$\Gamma$, and let $y\in\bar S$ with $\bar h(y)<(n+1)\delta$.  Choose $i$ with
nonzero $i$th coordinate.  Every summand in a sequence of sums converging to
$y-g_i$ in \eqref{eq:descent} eventually has height less than $n\delta$, so
it lies in the finite set $\Gamma\cap\{\bar h<n\delta\}$.  Hence
$y-g_i\in\Gamma$ and $y\in\Gamma$.
\end{proof}

When the face is a single exposed point $v$, translated to the origin, we
have $V=\{0\}$, $N=\R^d$ and $\pi$ is the identity, so
\cref{lem:normal-basis} and \cref{thm:free-monoid} describe the
whole support; this is \cref{thm:boundary-intro}\textup{(c)}.

\subsection{Positive-Mass Facets}

\Cref{thm:massive-face} produces an exposed face of positive mass,
possibly of small dimension.  We now enlarge it to a facet and quantize the
normal coordinate there.

\begin{theorem}[Face Flag]\label{thm:face-flag}
For the face $F$ fixed above, $\aff F=\aff\supp(\mu|_F)$, and there is a flag of exposed
positive-mass faces $F=F_0\subset F_1\subset\dots\subset F_{k-1}$ with
$\dim F_j=r+j$.  In particular $F_{k-1}$ is a facet and
$\mu|_{F_{k-1}}$ spans its supporting hyperplane affinely.
\end{theorem}

\begin{proof}
Let $(g_j^*)$ be the dual basis of $N^*$ and
$\ell_j=g_{j+1}^*+\dots+g_k^*$.  By \eqref{eq:free-monoid}, every coordinate
of every point of $\bar S$ is a nonnegative integer, so
$\pi(Q)\subset\{\sum_ic_ig_i:c_i\ge0\}$.  Since $h=\bar h\circ\pi$ and
$\bar h(g_i)>0$, $\pi$ vanishes on $F$; thus $F\subset V$, which gives
\begin{equation*}
 \aff F=\aff\supp(\mu|_F),\qquad F=Q\cap\{\ell_0\circ\pi=0\}=F_0.
\end{equation*}
Moreover
$\ell_j\circ\pi\ge0$ on $Q$, so $F_j=Q\cap\{\ell_j\circ\pi=0\}$ is exposed
and contains $F$.
It also contains the sets $\supp\mu\cap\pi^{-1}(g_i)$, $i\le j$, which have
positive mass; their affine span together with the face has direction space
$V+\Span(g_1,\dots,g_j)$, of dimension $r+j$.  The $k-j$ supporting hyperplanes give the matching
upper bound.
\end{proof}

At a facet the normal quotient is one-dimensional, and the monoid becomes
an arithmetic progression.

\begin{theorem}[Facet Quantization]\label{thm:facet}
Let $F=Q\cap\{h=0\}$, $h=\langle a,\cdot\rangle-\alpha\ge0$ on $Q$, be a
facet with $\mu(F)>0$, and
let $g=\langle a,b\rangle$.  Then $g>0$, $\supp h_\#\mu\subset g\Nzero$, and
$h_\#\mu$ has atoms at $0$ and $g$.
\end{theorem}

\begin{proof}
By \cref{thm:face-flag}, $\supp(\mu|_F)$ spans the hyperplane
$\{h=0\}$ affinely, so $k=1$.
\Cref{lem:normal-basis} gives $g>0$, and gives an atom at $g$ because
$\bar h(g_1)=\langle a,b\rangle$.  Scale $h$ so that $g=1$ and let
$S=\supp h_\#\mu$.  Every $s\in S\setminus\{0\}$ is the height of a tuple
consisting of one support point at height near $s$ and $d$ affinely
independent face points; such a tuple is affinely independent.  Projecting
\eqref{eq:support-identity} by the sum of heights gives
$S\setminus\{0\}\subset1+\overline{S^{\oplus(d+1)}}$.  By induction,
$S\cap[0,n)\subset\{0,\dots,n-1\}$.  Indeed, if $s\in S\cap(n,n+1)$, then
$s-1$ is a limit of sums of $d+1$ points of $S$ below $n$.  These points are
integers by the induction hypothesis, so the sums take finitely many
integer values and cannot converge to $s-1\in(n-1,n)$.
\end{proof}

\begin{proof}[Proof of \cref{thm:boundary-intro}]
By \cref{thm:constant-intro}, a non-Gaussian solution has $b\ne0$;
conversely a Gaussian solution has $b=0$.  Part~(a) follows from
\cref{cor:face-null}(i), \cref{thm:face-flag} and
\cref{thm:facet}.  Part~(b) is \cref{lem:normal-basis} with
\cref{thm:free-monoid}.  Part~(c) is the case $F=\{v\}$ of part~(b),
in which $N=\R^d$ and $\pi$ is the identity after $v$ is translated to the
origin; the linear map sending $g_i$ to the $i$th basis vector then puts the
support in $\Nzero^d$.  The final assertions are
\cref{cor:face-null}.
\end{proof}

\section{Cumulant Rigidity and Infinitely Divisible Solutions}
\label{sec:cumulant}

This section proves \cref{thm:cum-intro} and
\cref{cor:id-intro}.  The boundary argument used positivity
globally through supports; we now use it locally through cumulants.  Under
\eqref{eq:E}, a natural trace of the
tensor $\mathscr A_\theta$ defined below vanishes identically; in the
geometric reading of \cref{sec:kahler} it is a Ricci curvature.  A
one-sided sign condition therefore forces the tensor itself to vanish.

Put $H=\nabla^2K$ and let $(H^{ab})=H^{-1}$.  Subscripts on $K$ denote
partial derivatives, so that $K_{ab}=H_{ab}$, while $K_{abc}$ and $K_{abcd}$
are the third and fourth cumulant tensors of $P_\theta$; repeated indices
are summed.  On symmetric matrices define
\begin{align*}
 \mathscr A_\theta(S,T)&=A_{ab,cd}S^{ab}T^{cd},&
 A_{ab,cd}&=K_{abcd}-K_{abp}H^{pq}K_{qcd},\\
 \mathscr W_\theta(S,T)&=W_{ab,cd}S^{ab}T^{cd},&
 W_{ab,cd}&=H_{ac}H_{bd}+H_{ad}H_{bc},\\
 B_\theta(u,v)&=\mathscr A_\theta(u\odot v,u\odot v),
\end{align*}
where $u\odot v=(uv^{\mathsf T}+vu^{\mathsf T})/2$ is the symmetrized
product of $u$ and $v$.  Thus, for $X\sim P_\theta$ and
$Y_u=\langle u,X\rangle$,
\[
 B_\theta(u,v)=\operatorname{cum}(Y_u,Y_v,Y_u,Y_v)
 -c^{\mathsf T}\Cov(X)^{-1}c,\qquad
 c=\operatorname{cum}(Y_u,Y_v,X).
\]
We next recall the J{\o}rgensen set of $\mu$, the set of orders $r$ for
which the $r$th convolution power of $\mu$ exists \cite{Jorgensen1987},
\cite[Definition~3.1]{KokonendjiSeshadri1996}:
\[
 \Lambda(\mu)=\{r>0:Z^r=L_{\mu_r}\text{ on }\Theta
 \text{ for a positive Radon measure }\mu_r\}.
\]
It contains the positive integers, and it is closed under addition because
convolution multiplies Laplace transforms.  It is also closed in
$(0,\infty)$.  Indeed, let $r_n\in\Lambda(\mu)$ converge to $r>0$, fix
$\theta_0\in\Theta$, and let $Q_n$ be the probability measure
$Z(\theta_0)^{-r_n}\ee^{\langle\theta_0,x\rangle}\mu_{r_n}(\dd x)$.  Its
Laplace transform is $m^{r_n}$, where $m(s)=Z(\theta_0+s)/Z(\theta_0)$ on
$U=\Theta-\theta_0$.  These transforms are bounded near the origin uniformly
in $n$, so the sequence $(Q_n)$ is tight.  For $s\in U$ choose $p>1$ with
$ps\in U$; the bound $m(ps)^{r_n}$ makes $\ee^{\langle s,x\rangle}$ uniformly
integrable, so a weak limit $Q$ of a subsequence has Laplace transform $m^r$
on $U$, and $\mu_r(\dd x)=Z(\theta_0)^r\ee^{-\langle\theta_0,x\rangle}Q(\dd x)$
has Laplace transform $Z^r$ on $\Theta$.  Albeverio, Gottschalk and Wu prove
the same properties for the analogous set defined through characteristic
functions of laws on $\R$ \cite[Proposition~2.1]{AlbeverioGottschalkWu1998}.
The set $\Lambda(\mu)$ equals $(0,\infty)$ precisely when the laws $P_\theta$
are infinitely divisible.  Otherwise $\inf\Lambda(\mu)>0$, because a closed
additive semigroup that contains $1$ and arbitrarily small positive numbers
is all of $(0,\infty)$ \cite[Proposition~2.2]{AlbeverioGottschalkWu1998}.
On $M=\nabla K(\Theta)$ write
$V_\mu(m)=H((\nabla K)^{-1}(m))$ for the variance function.  The following
lemma collects the identities on which the section rests; for quadratic
forms on symmetric matrices, $\succeq0$ means positive semidefinite.

\begin{lemma}[Trace, Variance, and Power Identities]\label{lem:cum-master}
For an $H$-orthonormal basis $(e_k)$ and $m=\nabla K(\theta)$,
\begin{align}
 \sum_kB_\theta(u,e_k)&=\partial_u^2\log\det H,\label{eq:cum-trace}\\
 A_{ab,cd}&=H_{ak}H_{bl}\,
 \frac{\partial^2(V_\mu)_{cd}}{\partial m_k\partial m_l}.
 \label{eq:cum-variance}
\end{align}
Moreover, $B_\theta\equiv0$ implies $\mathscr A_\theta=0$, and
$\mathscr A_\theta+r\mathscr W_\theta\succeq0$ whenever
$r\in\Lambda(\mu)$.
\end{lemma}

\begin{proof}
Jacobi's formula gives
\begin{equation*}
 \partial_a\partial_c\log\det H=H^{bd}K_{abcd}-H^{bd}H^{pq}K_{abp}K_{cdq},
\end{equation*}
and $\sum_ke_k^be_k^d=H^{bd}$; this is \eqref{eq:cum-trace}.  To see that
$B_\theta\equiv0$ implies $\mathscr A_\theta=0$, note that $A_{ab,cd}$ is
symmetric in $(a,b)$, in
$(c,d)$ and under exchange of the two pairs.  If
$B_\theta(u,v)=A_{ab,cd}u^av^bu^cv^d$ vanishes identically, polarizing
separately in $u$ and in $v$ gives $A_{ab,cd}=-A_{cb,ad}$.  Thus the tensor
is fixed by the transposition $(1\,2)$ of its index positions and negated by
$(1\,3)$.  Since $(2\,3)=(1\,2)(1\,3)(1\,2)$ and $(1\,2)=(1\,3)(2\,3)(1\,3)$,
it is negated by $(2\,3)$ and then by $(1\,2)$, so it is zero.
Differentiating $V_\mu(\nabla K)=H$ once gives
$\partial_{m_k}(V_\mu)_{cd}=K_{cdp}H^{pk}$, and differentiating again gives
\eqref{eq:cum-variance}.  For the last assertion, fix $r\in\Lambda(\mu)$
and let $X$ have the law
$\exp\{\langle\theta,x\rangle-rK(\theta)\}\,\mu_r(\dd x)$, the $r$th
convolution power of $P_\theta$, and put $\widetilde X=X-\mathbb EX$.
Then $r(\mathscr A_\theta+r\mathscr W_\theta)$ is the covariance form of
the residual left after regressing the centred quadratic statistic
$\widetilde X\widetilde X^{\mathsf T}-rH$ on $\widetilde X$, so it is
positive semidefinite; \cref{app:residual} gives the computation.
\end{proof}

Under \eqref{eq:E} the right side of \eqref{eq:cum-trace} vanishes.  This
turns a sign condition on $B_\theta$ into rigidity.

\begin{theorem}[Localized Fourth-Cumulant Criterion]\label{thm:cum}
Assume \eqref{eq:E}.  Then $\mu$ is Poisson--Gaussian if and only if
$B_\theta(u,v)\ge0$ for all $u,v$ throughout some nonempty open subset of
$\Theta$.  These conditions are also equivalent to $\mathscr A\equiv0$
and to $V_\mu$ being affine.  If
$\mu$ is not Poisson--Gaussian, the set
\[
 \{\theta\in\Theta:B_\theta(u,v)<0\text{ for some }u,v\}
\]
is open and dense.
\end{theorem}

\begin{proof}
The Poisson--Gaussian variance function is affine.  Since $H$ is invertible,
\eqref{eq:cum-variance} identifies affine $V_\mu$ with
$\mathscr A=0$, which gives $B=0$; the converse classification is the
affine case of Letac's theorem \cite{Letac1989}, or, in Hessian-geometric
terms, \cite[Theorem~4.1]{FuruhataKurose2013}.  Now suppose that $B\ge0$
on the open set.  Equations \eqref{eq:E} and \eqref{eq:cum-trace} give
$0=\sum_kB_\theta(u,e_k)$.  Nonnegativity makes every summand zero; varying
the orthonormal frame gives $B_\theta\equiv0$, and
\cref{lem:cum-master} then gives
$\mathscr A=0$ on the open set, hence everywhere by analyticity.  If the
family is not Poisson--Gaussian, the analytic tensor $\mathscr A$ is nonzero
on an open dense set; at every such point \eqref{eq:cum-trace} forces $B$ to
take a negative value.  A direct proof of the affine classification in the
present normalization is given in \cref{prop:app-affine}.
\end{proof}

By the cumulant expression of $B_\theta$, \cref{thm:cum} contains
\cref{thm:cum-intro}.  Infinitely divisible families satisfy the sign
condition, which proves \cref{cor:id-intro}.

\begin{proof}[Proof of \cref{cor:id-intro}]
Infinite divisibility gives $\Lambda(\mu)=(0,\infty)$.  Letting $r\downarrow0$
in \cref{lem:cum-master} gives $\mathscr A\succeq0$, hence $B\ge0$ on
$\Theta$; \cref{thm:cum} applies.
\end{proof}

For infinitely divisible laws in dimension one, the condition $B\ge0$ reads
$(K''')^2\le K''K''''$; this is the inequality between squared skewness and
excess kurtosis of Rohatgi and Sz\'ekely \cite{RohatgiSzekely1989}.  After
standardizing $H(\theta)=I$, evaluating
$\mathscr A_\theta\succeq0$ at the identity recovers the multivariate
inequality of M\'ori, Rohatgi and Sz\'ekely
\cite[Theorem~2]{MoriRohatgiSzekely1993}; related multivariate inequalities
are studied by Ogasawara~\cite{Ogasawara2017}.  Equation \eqref{eq:E} is
essential in \cref{thm:cum}: \cref{ex:app-Bernoulli} gives a
family, not satisfying \eqref{eq:E}, for which $B_\theta\ge0$ at every
$\theta$ although it is not infinitely divisible.

Finally, \cref{cor:id-intro} constrains the solutions that are not
Poisson--Gaussian.  Such a solution is not infinitely divisible, so
$\Lambda(\mu)\ne(0,\infty)$, and the semigroup property recalled above gives
\begin{equation}\label{eq:jorgensen-gap}
 \inf\Lambda(\mu)>0 .
\end{equation}
In other words, a solution that is not Poisson--Gaussian has no convolution
powers of arbitrarily small order.

\section{Projectively Induced Ricci-Flat K\"ahler Tubes}
\label{sec:kahler}

This section proves the geometric statements of
\cref{sec:intro}.  The tensor $\mathscr A$ and the J{\o}rgensen set
have direct geometric meanings for the K\"ahler metric $\omega_K$ on the
tube over $\Theta$: the former is, up to normalization and the order of the
indices, its curvature tensor, while the latter records the scales at which
the metric is projectively induced.  This is the tube analogue of
the description, due to Loi and Zedda, of the projectively induced multiples
of the Bergman metric of a bounded symmetric domain through its Wallach set
\cite[Theorem~2]{LoiZedda2011}.

Let $\Theta\subset\R^d$ be convex and open and let $K\in C^\infty(\Theta)$
be strictly convex.  On the tube $T\Theta=\Theta+i\R^d$ consider the
K\"ahler form $\omega_K=i\partial\bar\partial K(\operatorname{Re}z)$.  It is
half of the K\"ahler form that Dombrowski's construction associates, on the
tangent bundle of $\Theta$, with the Hessian metric $\nabla^2K$ and the flat
connection \cite{Dombrowski1962}, \cite[Chapter~2]{Shima2007}.  Its Ricci
form is
$-i\partial\bar\partial\log\det\nabla^2K$.  Hence $\omega_K$ is Ricci-flat if
and only if $\log\det\nabla^2K$ is affine, i.e.\ \eqref{eq:E} holds.  For
exponential families, this equivalence appears in
\cite[Proposition~14]{ZhangKhan2020}.  If
$K=\log\int\ee^{\langle\theta,x\rangle}\mu(\dd x)$, then with the
normalization $\omega_{\mathrm{FS}}=i\partial\bar\partial\log\|v\|^2$, the
coherent-state map $z\mapsto[x\mapsto\ee^{\langle z,x\rangle/2}]$ is a
K\"ahler immersion of $(T\Theta,\omega_K)$ into the projective Hilbert space
$\mathbb P(L^2(\mu))$.  Special cases of this construction appear
in~\cite[Proposition~7.2 and Sections~9--10]{Molitor2013}.

Write $g_{i\bar j}=K_{ij}/4$ for the metric of $\omega_K$.  With the
curvature convention
\begin{equation*}
 R_{i\bar jk\bar l}=-\partial_i\partial_{\bar j}g_{k\bar l}
 +g^{p\bar q}(\partial_i g_{k\bar q})(\partial_{\bar j}g_{p\bar l}),
\end{equation*}
a direct calculation gives
\begin{equation}\label{eq:kahler-curvature}
 R_{i\bar jk\bar l}
 =-\frac1{16}\bigl(K_{ijkl}-K_{ikp}H^{pq}K_{qjl}\bigr),
 \qquad
 B_\theta(u,v)=-16R(u,\bar u,v,\bar v)\quad(u,v\in\R^d).
\end{equation}
Up to index order and normalization, this is the curvature tensor of Hessian
geometry; see Shima~\cite[Chapter~3]{Shima2007}.  The trace identity
\eqref{eq:cum-trace} is therefore the usual statement that Ricci curvature is
the trace of holomorphic bisectional curvature, and
\cref{thm:cum} has the following geometric reading: a
Ricci-flat tube with $K=\log L_\mu$ whose holomorphic bisectional curvature
is nonpositive on the directions spanned by real tangent vectors is flat.

The coherent-state construction has a converse: every projective K\"ahler
immersion of a tube arises from a positive measure.

\begin{proposition}[Projective Reconstruction]\label{prop:reconstruction}
Let
\[
 F:(T\Theta,\omega_K)\longrightarrow
 (\mathbb P(\mathcal H),\omega_{\mathrm{FS}})
\]
be a K\"ahler immersion into a finite- or infinite-dimensional projective
Hilbert space.  Then there is a positive Radon measure $\mu$ on $\R^d$ with full
affine support such that $\ee^{K(\theta)}=\int\ee^{\langle\theta,x\rangle}
\mu(\dd x)$ on $\Theta$.  After restricting the target to the closed
projective span of the image, $F$ is projectively unitarily congruent to the
coherent-state map of $\mu$.
\end{proposition}

\begin{proof}
The tube is convex, hence Stein and contractible, so $F$ has a
nowhere-vanishing holomorphic lift $f_0$.  Writing
$\log\|f_0\|^2-K(\operatorname{Re}z)=\operatorname{Re}h$ with $h$
holomorphic, the lift $f=\ee^{-h/2}f_0$ satisfies
$\|f(z)\|^2=\ee^{K(\operatorname{Re}z)}$.  Its Gram kernel
$\langle f(z),f(w)\rangle$ is positive definite, sesquiholomorphic and, by
polarization, invariant under simultaneous imaginary translations.  The
scalar case of \cite[Theorem~III.5]{Neeb1998} writes it as
$\int\ee^{\langle z+\bar w,\xi\rangle}\nu(\dd\xi)$ with $\nu\ge0$; on real
points this is Devinatz's theorem~\cite[Theorem~1]{Devinatz1955}, which is
stated for planar domains with a remark that the proof extends to higher
dimensions.  The image $\mu$ of
$\nu$ under $\xi\mapsto2\xi$ has full affine support because $\nabla^2K>0$.
The coherent vectors are total in $L^2(\mu)$
\cite[Lemma~III.8(b)]{Neeb1998}, and Calabi's rigidity~\cite{Calabi1953}
gives the congruence.
\end{proof}

By analytic continuation, the same holds for immersions defined near a
single point and for every rescaling of the metric.

\begin{corollary}[Local Projective Consequences]
\label{cor:local-reconstruction}
Let $K$ be real analytic.  If $c\,\omega_K$, $c>0$, admits a K\"ahler
immersion near one point into a finite- or infinite-dimensional projective
Hilbert space, then
\[
 \ee^{cK(\theta)}
 =\int_{\R^d}\ee^{\langle\theta,x\rangle}\,\nu(\dd x)
\]
on $\Theta$ for a positive Radon measure $\nu$ with full affine support.
If moreover $K=\log L_\mu$ satisfies \eqref{eq:E} and $\mu$ is not
Poisson--Gaussian, then no $r\,\omega_K$ with $0<r<\inf\Lambda(\mu)$ admits
even a local projective K\"ahler immersion.
\end{corollary}

\begin{proof}
The tube is simply connected, so the local immersion extends to the whole
tube by Calabi's criterion \cite[Theorems~11--12]{Calabi1953}; for path
extension into projective Hilbert space see \cite[Theorem~1.1]{HuangLi2025}.
\Cref{prop:reconstruction}, applied to $cK$, then gives
$\nu$.  In the second assertion, a local immersion of $r\,\omega_K$ would
give $r\in\Lambda(\mu)$, and $\inf\Lambda(\mu)>0$ by
\eqref{eq:jorgensen-gap}.
\end{proof}

When $b=0$, the holomorphic volume form $\dd z_1\wedge\dots\wedge\dd z_d$
is parallel and the vertical fibres $\{\theta\}+i\R^d$ are special
Lagrangian.  After quotienting the imaginary directions by a lattice, they
become the torus fibres of the semi-flat Calabi--Yau model over the regular
part of a special Lagrangian torus fibration~\cite{LeungYauZaslow2000}.
Conversely, if every fibre is special Lagrangian for a parallel holomorphic
volume form $f\,\dd z_1\wedge\dots\wedge\dd z_d$, then
$|f|^2$ is a constant multiple of $\det\nabla^2K(\operatorname{Re}z)$ and
the phase of $f$ is constant along each fibre.  Hence $f$ does not depend on
$\operatorname{Im}z$, the Cauchy--Riemann equations make $f$ constant, and
$\det\nabla^2K$ is constant.  Accordingly we call $\omega_K$
\emph{semi-flat} when $\det\nabla^2K$ is constant.  In this case a
projective K\"ahler immersion at a single scale forces flatness.

\begin{theorem}[Semi-Flat Rigidity]\label{thm:semiflat}
Let $K$ be strictly convex on the convex open set $\Theta$ with
$\det\nabla^2K$ constant.  Suppose that for some $c>0$ the metric
$c\,\omega_K$ admits a K\"ahler immersion, possibly defined only near one
point, into a finite- or infinite-dimensional projective Hilbert space.
Then $K$ is a quadratic polynomial and $\omega_K$ is flat.  Conversely,
every positive-definite quadratic $K$ has a Gaussian coherent-state
immersion at every scale.
\end{theorem}

\begin{proof}
The equation $\det\nabla^2K=\mathrm{const}$ is analytic, and it is elliptic
along $K$: its linearization $\cof(\nabla^2K)^{ij}\partial_i\partial_j$ has
a positive-definite coefficient matrix because $\nabla^2K$ is positive
definite.  Morrey's
theorem \cite[p.~198]{Morrey1958}, as quoted in the proof of
\cref{prop:exp-moments}, now with the single unknown $u=K$, $s_1=0$ and
$t_1=2$, shows that the smooth function $K$ is real analytic.
\Cref{cor:local-reconstruction} gives $\mu$ with
$cK=\log L_\mu$ on $\Theta$.  Let $\kappa_0$ be the constant value of
$\det\nabla^2K$.  The analytic
function $\det\nabla^2\log L_\mu-c^d\kappa_0$ vanishes on $\Theta$, hence on
the whole interior of the natural parameter domain of $\mu$.
\Cref{thm:constant-intro} makes $\mu$ Gaussian and $K$ quadratic.
Conversely, a positive-definite quadratic $K$ is the logarithm of a Gaussian
Laplace transform, whose coherent-state map is defined at every scale.
\end{proof}

For $b\ne0$ the conclusion of \cref{thm:semiflat} fails as stated:
the Poisson tube $K=\ee^{\theta}$ is, up to a constant, the tube of a
cumulant function, hence projectively induced by its coherent-state map,
and $K$ is not quadratic; whether $\omega_K$ must still be flat is the
conjecture of Loi, Salis and Zuddas recalled in \cref{rem:context}.

\begin{remark}[Context]\label{rem:context}
Loi, Salis and Zuddas conjectured that every Ricci-flat projectively induced
K\"ahler metric is flat \cite[Conjecture~1]{LoiSalisZuddas2018}.
Previously settled cases include radial and stable-projectively-induced
metrics \cite{LoiSalisZuddas2018}, the Taub--NUT metrics
\cite{LoiZeddaZuddas2012,Huang2026TaubNut},
Calabi's Ricci-flat metrics on
canonical bundles \cite{LoiZeddaZuddas2021}, some of Stenzel's metrics
\cite{Zedda2021}, cones over regular complete Sasakian manifolds
\cite{MariniTardiniZedda2024} and radial metrics at a single scale
\cite{ArezzoLoiPlaciniZedda2026}; none of these classes contains a non-flat
semi-flat tube.  Finite-dimensional targets are already excluded: compact
Ricci-flat submanifolds by Hulin~\cite{Hulin2000}, and the remaining ones by
Arezzo, Li and Loi~\cite[Corollary~1.7]{ArezzoLiLoi2025}, building on
Hulin~\cite{Hulin1996}.  Complete semi-flat tubes
are flat by Calabi's affine-sphere theorem
\cite[Corollary to Theorem~1]{Calabi1958}, \cite[Theorem~2]{Loftin2008}; hence \cref{thm:semiflat} concerns the
incomplete case, permits a $\mathbb{CP}^{\infty}$ target, and requires only a
germ.  A semi-flat quotient $T\Theta/i\mathcal P$ by a nonzero
discrete group $\mathcal P\subset\R^d$ of imaginary translations admits no
global projective immersion: a lift would be congruent to the coherent-state map of
a Gaussian measure, while periodicity would make
$\ee^{i\langle\tau,x\rangle/2}$ constant almost everywhere for every
$\tau\in\mathcal P$; compare \cite[Exercise~3.6.6]{LoiZedda2018}.
\end{remark}

\section{Further Remarks}
\label{sec:remarks}

We close by comparing our results with two earlier arguments in the
literature.

\begin{remark}[Kokonendji's 2003 Argument]\label{rem:2003}
\cite[Theorem~3.2]{Kokonendji2003} states that the solutions of \eqref{eq:E}
are exactly the Poisson--Gaussian families, with no assumption on $D$ or on
infinite divisibility.  The proof applies a rigidity theorem for entire
solutions of $\det f''=\mathrm{const}$ to a function $f$ with
$\det f''=\det(\ee^{-\langle b,\theta\rangle/d}\nabla^2K)$.  For $b=0$ one can
take $f=K$, and the only gap, that $\Theta$ need not be $\R^d$, is closed by
\cref{thm:regularity-intro,thm:constant-intro}.  For $b\ne0$
the matrix $\ee^{-\langle b,\theta\rangle/d}\nabla^2K$ is not a Hessian in
general, and the step would force it to be constant, which fails for two
independent unit Poisson laws, where it equals
$\diag(\ee^{(\theta_1-\theta_2)/2},\ee^{(\theta_2-\theta_1)/2})$.  Thus the
argument does not cover $b\ne0$, the remaining part of the case
$\beta=0$ of Casalis's question \cite[p.~1831]{Casalis1996}; the case
$\beta\ne0$ is also recorded as open in \cite[p.~359]{Kokonendji2003}.  None of
our results uses \cite[Theorem~3.2]{Kokonendji2003}.
\end{remark}

\begin{remark}[The Kokonendji--Masmoudi Lemma]\label{rem:KM-lemma}
The coordinate description of the L\'evy measure in
\cite[Lemma~7]{KokonendjiMasmoudi2006} is false as stated.  Write $\Sigma$
for the Gaussian covariance and $\Pi$ for the L\'evy measure of an
infinitely divisible law.  For
$X=(P,P+G)$, where $P$ and $G$ are independent, $P$ is Poisson and $G$ is
standard Gaussian, one has $\Sigma=\diag(0,1)$ while $\Pi$ is supported at
$(1,1)$.  For $X=(P_1+P_2,P_2)$, where $P_1,P_2$ are independent Poisson
random variables, one has $\Sigma=0$ while $\Pi$ has the nonorthogonal atoms
$(1,0)$ and $(1,1)$.  Both laws are infinitely divisible solutions of
\eqref{eq:E}, with $b=(1,1)$ and $b=(2,1)$, and both satisfy the hypothesis
of the lemma, condition~(8) of \cite{KokonendjiMasmoudi2006}, whose
coefficients are only required to be nonnegative.  The lemma is reused in the
proof of
\cite[Theorem~4.1]{NisaKokonendjiSaefuddin2015}.  The proof of
\cref{cor:id-intro} in \cref{sec:cumulant} does not use that
lemma.
\end{remark}

\section*{Acknowledgments and Tool Disclosure}

\noindent\textit{AI-use disclosure.}
The authors used Anthropic Claude Code and OpenAI Codex as interactive
research and writing tools.  These tools assisted with exploratory
discussion, testing and refinement of ideas, literature and source
organization, code development and verification, mathematical error
checking, and editorial revision.  The authors take full responsibility
for the content presented.

\appendix
\crefalias{section}{appendix}
\crefalias{subsection}{subappendix}
\section{Complements to the Cumulant Criterion}
\label{app:cumulant}

We retain the notation of \cref{sec:cumulant}.  This appendix proves
the last assertion of \cref{lem:cum-master}, gives a direct proof of
the classification step in \cref{thm:cum}, and gives an example
showing that \eqref{eq:E} cannot be dropped from that theorem.

\subsection{The Residual-Covariance Identity}
\label{app:residual}

We prove the regression identity stated in the proof of
\cref{lem:cum-master}.

\begin{proof}[Proof of the last assertion of \cref{lem:cum-master}]
Choose the measure $\mu_r$ represented by $r\in\Lambda(\mu)$, and let $X$
have law
\[
 P_{\theta,r}(\dd x)
 =\exp\{\langle\theta,x\rangle-rK(\theta)\}\,\mu_r(\dd x).
\]
The cumulant function of $X$ near the origin is
$rK(\theta+\zeta)-rK(\theta)$.  If
$\widetilde X=X-\mathbb EX$, its centred moments through order four are
\begin{align}
 \mathbb E[\widetilde X_a\widetilde X_b]
 &=rK_{ab},\label{eq:app-moment2}\\
 \mathbb E[\widetilde X_a\widetilde X_b\widetilde X_c]
 &=rK_{abc},\label{eq:app-moment3}\\
 \mathbb E[\widetilde X_a\widetilde X_b
             \widetilde X_c\widetilde X_d]
 &=rK_{abcd}
  +r^2(K_{ab}K_{cd}+K_{ac}K_{bd}+K_{ad}K_{bc}).
 \label{eq:app-moment4}
\end{align}
Put
\[
 \alpha_{ab}^{\ q}=K_{abp}H^{pq},
 \qquad
 \mathcal R_{ab}=\widetilde X_a\widetilde X_b-rH_{ab}
              -\alpha_{ab}^{\ q}\widetilde X_q.
\]
The choice of $\alpha$ is the least-squares coefficient, since
\eqref{eq:app-moment2}--\eqref{eq:app-moment3} give
\[
 \mathbb E[(\widetilde X_a\widetilde X_b-rK_{ab})\widetilde X_c]
 =rK_{abc}
 =\alpha_{ab}^{\ q}\mathbb E[\widetilde X_q\widetilde X_c].
\]
Expanding the covariance gives
\begin{align*}
 \mathbb E[\mathcal R_{ab}\mathcal R_{cd}]
={}&rK_{abcd}
  +r^2(K_{ab}K_{cd}+K_{ac}K_{bd}+K_{ad}K_{bc})\\
 &-r^2K_{ab}K_{cd}-r^2K_{ab}K_{cd}+r^2K_{ab}K_{cd}\\
 &-r\alpha_{cd}^{\ q}K_{abq}
  -r\alpha_{ab}^{\ q}K_{qcd}
  +r\alpha_{ab}^{\ p}\alpha_{cd}^{\ q}K_{pq}.
\end{align*}
The coefficient of $r^2K_{ab}K_{cd}$ is $1-1-1+1=0$.  Each of the first
two contractions in the last line equals
$K_{abp}H^{pq}K_{qcd}$, and the final contraction restores one copy of the
same quantity.  Consequently
\begin{equation}\label{eq:app-residual-covariance}
 \mathbb E[\mathcal R_{ab}\mathcal R_{cd}]
 =r\bigl(A_{ab,cd}+rW_{ab,cd}\bigr).
\end{equation}
For a symmetric matrix $S$, contracting
\eqref{eq:app-residual-covariance} with $S^{ab}S^{cd}$ gives the variance of
$\mathcal R_{ab}S^{ab}$, which is nonnegative.  Division by $r>0$ gives
$\mathscr A_\theta(S,S)+r\mathscr W_\theta(S,S)\ge0$.
\end{proof}

\subsection{The Affine-Variance Case}
\label{app:affine}

We next give a direct proof of the classification step in
\cref{thm:cum}, in the present normalization.

\begin{proposition}[The Affine-Variance Case]\label{prop:app-affine}
If $\mathscr A\equiv0$ on $\Theta$, equivalently if the variance function
$V_\mu$ is affine, then, after an invertible affine change of the
statistic, $K$ takes the form \eqref{eq:PG} up to an additive constant; that
is, $\mu$ is Poisson--Gaussian.
\end{proposition}

\begin{proof}
The two hypotheses are equivalent by \eqref{eq:cum-variance}.  Suppose
$\mathscr A=0$ and set
$\mathcal C_c=H^{-1}\partial_cH$.  Direct differentiation gives
\[
 (H\partial_d\mathcal C_c)_{ab}=A_{ad,bc}=0.
\]
The matrices $\mathcal C_c$ are therefore constant; equality of mixed
derivatives makes them commute.  Translate the parameter so that a chosen
point of $\Theta$ is the origin, and make a linear change of the statistic after which $H(0)$
is the identity matrix.
The $\mathcal C_c$ are then commuting symmetric matrices and hence
simultaneously diagonal.  It follows that $H$ is diagonal.  Hessian
integrability,
$\partial_cH_{jj}=\partial_jH_{jc}=0$ for $c\ne j$, then gives
\[
 H_{jj}(\theta)=C_j\ee^{\beta_j\theta_j},\qquad C_j>0.
\]
After two integrations and normalization at the origin, the
moment-generating function factorizes.  Its one-dimensional factors are
Gaussian for $\beta_j=0$ and
Poisson with intensity $C_j/\beta_j^2$, followed by the affine change
$n\mapsto\alpha_j+\beta_jn$, for some $\alpha_j\in\R$, when
$\beta_j\ne0$.  Uniqueness of
moment-generating functions gives a product decomposition of the law at the
chosen point $\theta_0$.  Undoing the translation multiplies this law by
$\ee^{K(\theta_0)-\langle\theta_0,x\rangle}$, which keeps each Gaussian or
Poisson factor of the same kind up to a positive constant; hence $\mu$ is a
positive multiple of a product of the same kind.  The
same simultaneous diagonalization is used in
\cite[Theorem~4.1]{FuruhataKurose2013}, whose models
\cite[Example~2.3]{FuruhataKurose2013} are, in affine coordinates, the
Poisson--Gaussian potentials.
\end{proof}

\subsection{The Cumulant Inequality without the Determinant Equation}
\label{app:outside-class}

The following example lies outside \eqref{eq:E} and shows that, without the
determinant equation, the pointwise cumulant inequality does not imply
infinite divisibility.

\begin{example}[Convolution with a Bernoulli Law]
\label{ex:app-Bernoulli}
Let $\phi(x)=(2\pi)^{-1/2}\ee^{-x^2/2}$, and let $\nu$ be the infinitely
divisible law with no Gaussian part and with L\'evy measure
\[
 \Pi(\dd x)=\frac{\phi(x)}{x^2}\,\dd x
 \qquad(x\ne0).
\]
This is a L\'evy measure because
$\int(1\wedge x^2)\Pi(\dd x)<\infty$.  Choosing an arbitrary drift, its
cumulant function satisfies
\begin{equation}\label{eq:app-nu-variance}
 K_\nu''(\theta)=\int_{\R}\ee^{\theta x}\phi(x)\,\dd x
 =\ee^{\theta^2/2}.
\end{equation}
Let $B$ be independent of $\nu$ with
$\mathbb P(B=0)=\mathbb P(B=1)=1/2$, and set $\mu=\nu*\mathcal L(B)$.
Then $\mu$ is not infinitely divisible: its characteristic function vanishes
at $\pi$ because of the Bernoulli factor, whereas an infinitely divisible
characteristic function has no zeros.

Nevertheless $B_\theta(1,1)=K_\mu''''-(K_\mu''')^2/K_\mu''>0$ at every
$\theta$.  Indeed, put
\[
 v(\theta)=\frac{\ee^\theta}{(1+\ee^\theta)^2},
 \qquad
 w(\theta)=\ee^{-\theta^2/2}v(\theta).
\]
Then
\begin{equation}\label{eq:app-Kpp-w}
 K_\mu''(\theta)=\ee^{\theta^2/2}+v(\theta)
 =\ee^{\theta^2/2}(1+w(\theta)).
\end{equation}
Writing $p=\ee^\theta/(1+\ee^\theta)$ gives
\[
 v=p(1-p),
 \qquad v'=v(1-2p),
 \qquad v''=v(1-6p+6p^2).
\]
Hence $0<v\le1/4$, $|v'|\le v$ and $|v''|\le v$.  With
$x=|\theta|$,
\begin{align}
 |w'|
 &\le\frac14(1+x)\ee^{-x^2/2}\le\frac38,
 \label{eq:app-wprime-bound}\\
 |w''|
 &\le\frac14(1+2x+|x^2-1|)\ee^{-x^2/2}
 \le\frac34.\label{eq:app-wsecond-bound}
\end{align}
The last inequalities require no numerical approximation.  On
$0\le x\le1$,
$(1+x)\ee^{-x^2/2}\le3/2$ follows from
$1+x\le\frac32(1+x^2/2)$, while on $x\ge1$ the left-hand side is decreasing.
For \eqref{eq:app-wsecond-bound}, the numerator is
$2+2x-x^2\le3$ on $[0,1]$.  On $[1,\infty)$ it is $x^2+2x$, and
\[
 x^2\ee^{-x^2/2}\le\frac2e<1,
 \qquad
 2x\ee^{-x^2/2}\le2\ee^{-1/2}<2.
\]
Since $w\ge0$, \eqref{eq:app-Kpp-w} and
\eqref{eq:app-wprime-bound}--\eqref{eq:app-wsecond-bound} give
\begin{equation*}
 (\log K_\mu'')''
 =1+\frac{w''}{1+w}-\frac{(w')^2}{(1+w)^2}
 \ge1-|w''|-|w'|^2
 \ge1-\frac34-\frac9{64}
 =\frac7{64}>0.
\end{equation*}
Thus $K^{(4)}-(K^{(3)})^2/K''>0$, exactly the strict version of
$B_\theta(u,v)\ge0$ in dimension one.  At the same time the strict convexity
of $\log K_\mu''$ shows explicitly that this example does not satisfy
\eqref{eq:E}, which in dimension one would make $\log K''$ affine.

Taking the product with a nondegenerate Gaussian law on $\R^{d-1}$ gives an
example in every dimension.  Its tensor $\mathscr A$ has only the
one-dimensional component above, so
\begin{equation*}
 B_\theta(u,v)=A_{11,11}(u_1v_1)^2\ge0.
\end{equation*}
\end{example}

\section{Dimension Two}
\label{app:dimension-two}

This appendix sharpens the global boundary mechanism of
\cref{sec:boundary} in dimension two.  \Cref{sec:rank-two} proves that a
non-Gaussian solution either lives on a lattice orthant or has a convex
support that contains a line (\cref{thm:rank-two-intro}).
\Cref{sec:harmonic-boundary} introduces the axisymmetric harmonic
reduction and states rigidity in three regular regimes
(\cref{thm:harmonic-intro}), whose three parts are proved in
\cref{sec:finite-axial,sec:stieltjes-rows,sec:walls}.

\subsection{Lattice Orthant or Lineality}
\label{sec:rank-two}

In dimension two the facet $F$ of \cref{thm:boundary-intro}(a) is a segment, a half-line or a line.  We call the direction of $F$
\emph{tangential} and the transverse coordinate $h$ \emph{normal}.  Either
the tangential support is quantized as well, or the convex support contains
a whole line.

\begin{theorem}[Dimension Two]\label{thm:rank-two-intro}
Let $d=2$ and let $\mu$ be a non-Gaussian solution of \eqref{eq:E}.  Then
exactly one of the following holds.
\begin{enumerate}
 \item[\textup{(i)}] $\mu$ charges an exposed point of $Q$.  Hence, by
 \cref{thm:boundary-intro}\textup{(c)}, after an affine change of
 coordinates $\supp\mu\subset\Nzero^2$ and $\mu$ charges $(0,0)$, $(1,0)$
 and $(0,1)$.
 \item[\textup{(ii)}] $Q$ contains a line parallel to the facet of
 \cref{thm:boundary-intro}\textup{(a)} and is invariant under
 translations along it.
\end{enumerate}
\end{theorem}

Thus a non-Gaussian solution in the plane either is a bivariate
power-series family on the lattice orthant $\Nzero^2$, as for two
independent Poisson laws, or has a convex support that is invariant under a
line of translations, as for the product of a Poisson and a Gaussian law.
We call these the \emph{lattice-orthant branch} and the \emph{lineality
branch}; the lineality space of a closed convex set consists of the
directions of the lines that it contains
\cite[Section~8]{Rockafellar1970}.  Neither branch is classified here.

To prove \cref{thm:rank-two-intro}, we compare, in the measure
identity \eqref{eq:measure-identity}, the parts of $\mu$ at the successive
integer values of the normal coordinate.

In this subsection $d=2$ and $\mu$ is a non-Gaussian solution of
\eqref{eq:E}.  By \cref{thm:boundary-intro}(a), after an affine change of
coordinates $x=(n,y)$ and a rescaling of the normal coordinate, we may
assume that $Q\subset[0,\infty)\times\R$, that $F=Q\cap(\{0\}\times\R)$ is a
facet such that $\supp(\mu|_F)$ contains at least two points, that
$\supp\mu\subset\Nzero\times\R$, and that $b=(1,\beta)$.  Write
$\mu=\sum_{n\ge0}\delta_n\otimes\mu_n$ with positive measures $\mu_n$
on $\R$, so that $\mu_0\ne0$ and $\mu_1\ne0$.  We call $n$ the normal and
$y$ the tangential coordinate, $\mu_n$ the $n$th \emph{row} of $\mu$, and
$\mu_0$ the bottom row.

Let $\mathscr D_N$ be the restriction of $\nu_\mu$ to $\{N\}\times\R$,
regarded as a measure on $\R$.  Comparing normal levels in
\eqref{eq:measure-identity} gives, for $N\ge1$,
\begin{equation}\label{eq:level}
 \mathscr D_N=C\,\delta_\beta*\sum_{i+j+k=N-1}\mu_i*\mu_j*\mu_k ,
\end{equation}
summed over ordered triples.  For the points $(0,u)$, $(0,v)$ and $(N,y)$,
the determinant in \eqref{eq:nu-def} equals $\pm N(u-v)$.  Hence the
configurations with two points in row $0$ contribute
$N^2\nu_0*\mu_N\le\mathscr D_N$, where
$\nu_0(E)=\frac12\iint(u-v)^2\mathbf 1_{\{u+v\in E\}}\mu_0(\dd u)\mu_0(\dd
v)$.  For $N=1$ these are the only configurations, so
\begin{equation}\label{eq:leading-row}
 \nu_0*\mu_1=C\,\delta_\beta*\mu_0^{*3}.
\end{equation}

We begin with a dichotomy that uses only the convex geometry of $Q$.

\begin{lemma}[Pointed Support or Lineality]\label{lem:pointed}
Either $Q$ contains a line, which is then parallel to $F$ and along which $Q$
is invariant, or the following holds after a change of coordinates
$(n,y)\mapsto(n,\kappa n+\varepsilon y+c')$ with $\varepsilon\in\{\pm1\}$:
every row $\mu_n$ is supported in $[0,\infty)$, $0\in\supp\mu_0$, and the
origin is an exposed point of $Q$.
\end{lemma}

\begin{proof}
Since $n\ge0$ on $Q$, every line in $Q$ has direction $(0,\pm1)$, and a
closed convex set is invariant under translations in its lineality space
\cite[Section~8]{Rockafellar1970}.  If $Q$ contains
no line, then $F$ is a segment or a half-line.  Let $y_0$ be a finite
endpoint of $F$, and choose $\varepsilon$ with $F\subset\{\varepsilon(y-y_0)\ge0\}$.
Then
\[
 c_0:=-\inf\Bigl\{\frac{\varepsilon(y-y_0)}n:(n,y)\in Q,\ n\ge1\Bigr\}<\infty .
\]
Otherwise there would be points $(n_j,y_j)\in Q$ with
$\varepsilon(y_j-y_0)/n_j\to-\infty$.  Their directions from $(0,y_0)$
would converge to $(0,-\varepsilon)$, which would then be a recession
direction of $Q$ and would give points of $F$ with
$\varepsilon(y-y_0)<0$.  Put $y'=\kappa n+\varepsilon(y-y_0)$ with
$\kappa>c_0+1$.  Then $y'\ge0$ on row $0$, with equality only at $y_0$,
and $y'\ge(\kappa-c_0)n\ge1$ on the rows $n\ge1$.  The functional $n+y'$ is
nonnegative on $Q$ and vanishes only at the origin, which is therefore an
exposed point.  The origin lies in $Q$ and is a limit of convex combinations
of support points; since $n+y'\ge1$ off row $0$, row $0$ contains support
points with arbitrarily small $y'$, and $0\in\supp\mu_0$ because supports
are closed.
\end{proof}

In the second case of \cref{lem:pointed}, the vector $b$ becomes
$(1,\beta')$ with $\beta'=\kappa+\varepsilon\beta$.

We call the endpoint $0$ of the bottom row \emph{diffuse} if
$\mu_0\{0\}=0$.  In that case the supports of the rows drift linearly.

\begin{lemma}[Drift of a Diffuse Endpoint]\label{lem:drift}
In the second case of \cref{lem:pointed}, suppose that $\mu_0\{0\}=0$.
Then $\beta'\ge0$ and $\supp\mu_N\subset[N\beta',\infty)$ for all $N\ge0$.
\end{lemma}

\begin{proof}
Since $0\in\supp\mu_0$ and $\mu_0\{0\}=0$, distinct support points of
$\mu_0$ accumulate at $0$.  Hence $\inf\supp\nu_0=0$ and $\nu_0\{0\}=0$.
Comparing infima of supports in \eqref{eq:leading-row} gives
$\inf\supp\mu_1=\beta'$, and $\beta'\ge0$ because $\mu_1$ lives on
$[0,\infty)$.  If $\supp\mu_j\subset[j\beta',\infty)$ for $j<N$, then the
right side of \eqref{eq:level} is supported in $[N\beta',\infty)$.  Since
$N^2\nu_0*\mu_N\le\mathscr D_N$ and $\nu_0$ has mass arbitrarily close to
$0$, $\supp\mu_N\subset[N\beta',\infty)$.
\end{proof}

Let $G_N(s)=\int\ee^{sy}\mu_N(\dd y)$, so that $Z(t,s)=\sum_N\ee^{Nt}G_N(s)$.
Under the conclusion of \cref{lem:drift}, put
$\widetilde G_N(s)=\ee^{-N\beta's}G_N(s)$, the Laplace transform of $\mu_N$
translated by $-N\beta'$, and
$\widetilde Z(q,s)=\sum_N\ee^{Nq}\widetilde G_N(s)=Z(q-\beta's,s)$.  The
change of variables $(q,s)\mapsto(q-\beta's,s)$ has determinant one, so
\eqref{eq:E} becomes the normalized determinant equation
\begin{equation}\label{eq:normalized-determinant}
 \det D^2_{(q,s)}\log\widetilde Z=C\ee^q ,
\end{equation}
all rows of $\widetilde Z$ are supported in $[0,\infty)$, and the bottom row
is $\mu_0$.

The rescaling argument below needs a growth property of Laplace transforms
near a diffuse endpoint, which holds automatically.

\begin{lemma}[Automatic Endpoint Regularity]\label{lem:O-regular}
Let $\sigma$ be a positive measure on $[0,\infty)$ such that
$\int\ee^{-t_1y}\sigma(\dd y)<\infty$ for some $t_1>0$, with
$0\in\supp\sigma$ and $\sigma\{0\}=0$.  Put $L_\sigma(t)=\int\ee^{-ty}\sigma(\dd y)$ and
$\varphi=-\log L_\sigma$.  Then $\varphi(t)\to\infty$ and
$\varphi(t)/t\to0$ as $t\to\infty$.  Moreover, $\varphi$ is concave and
nondecreasing on $[t_1,\infty)$, and for every
compact $I\subset(0,\infty)$ there are $B,R_1<\infty$ with
$\varphi(Ru)\le B\varphi(R)$ for $u\in I$ and $R\ge R_1$.
\end{lemma}

\begin{proof}
Monotone convergence gives $L_\sigma(t)\to\sigma\{0\}=0$.  For every
$\delta>0$, $L_\sigma(t)\ge\sigma([0,\delta])\ee^{-\delta t}$, so
$\limsup\varphi(t)/t\le\delta$.  The function $\log L_\sigma$ is convex and
nonincreasing, so $\varphi$ is concave and nondecreasing on
$[t_1,\infty)$.  By concavity the chord slope
$(\varphi(t)-\varphi(t_1))/(t-t_1)$ is nonincreasing, so for $u\ge1$ and
$R\ge2t_1$ we get $\varphi(uR)\le|\varphi(t_1)|+2u(\varphi(R)+
|\varphi(t_1)|)$.  Monotonicity handles $u\le1$, and
$\varphi(R)\to\infty$ gives the last claim.
\end{proof}

We can now exclude a diffuse endpoint.  After rescaling, the logarithm of
$\widetilde Z$ converges to a convex function of one variable, whose
Monge--Amp\`ere measure vanishes, although \eqref{eq:normalized-determinant}
is invariant under the rescaling.

\begin{theorem}[No Diffuse Endpoint]\label{thm:no-diffuse}
In the second case of \cref{lem:pointed}, $\mu_0\{0\}>0$.
\end{theorem}

\begin{proof}
Suppose $\mu_0\{0\}=0$ and use the notation of
\eqref{eq:normalized-determinant},
writing $x=\ee^q$.  Since the rows are supported in $[0,\infty)$, each
$\widetilde G_N$ is nondecreasing.  Hence there are $x_0>0$ and $t_0>0$
with $M_0:=\widetilde Z(\log x_0,-t_0)<\infty$, and
$\widetilde G_N(s)\le M_0x_0^{-N}$ for $s\le-t_0$.  Let $\varphi$ be as in
\cref{lem:O-regular} with $\sigma=\mu_0$, and put $\Phi_R=\varphi(R)$
and $\varepsilon_R=\Phi_R^2/R^2\to0$.  Define
\[
 U_R(q,u)=\Phi_R^{-1}\log\widetilde Z(\log\varepsilon_R+q,-Ru).
\]
The Hessian of $U_R$ is obtained from that of $\log\widetilde Z$ by the
diagonal scaling $(1,-R)$ and the factor $\Phi_R^{-1}$.  By
\eqref{eq:normalized-determinant},
\begin{equation}\label{eq:scale-invariant}
 \det D^2U_R=\Phi_R^{-2}R^2C\varepsilon_R\ee^q=C\ee^q .
\end{equation}

Fix a rectangle $[q_-,q_+]\times[u_-,u_+]\subset\R\times(0,\infty)$.  For
large $R$, $\varepsilon_R\ee^{q_+}\le x_0/2$ and $Ru_-\ge t_0$, so $U_R$ is
finite there.  Moreover
\begin{equation*}
 -B\le-\varphi(Ru)/\Phi_R\le U_R\le\Phi_R^{-1}\log M_0,
\end{equation*}
using $\widetilde Z\ge\widetilde G_0(-Ru)=\ee^{-\varphi(Ru)}$ and
\cref{lem:O-regular}.  We claim that $\partial_qU_R\to0$ uniformly on
the rectangle.  Indeed,
\[
 \partial_qU_R=\Phi_R^{-1}\,
 \frac{\sum_NNx^N\widetilde G_N(-Ru)}{\widetilde Z(\log x,-Ru)},
 \qquad x=\varepsilon_R\ee^q .
\]
Let $\rho_R=\varepsilon_R\ee^{q_+}/x_0\to0$ and
$L_R=\log(1/\rho_R)\to\infty$, and split the sum at
$N_R=\lceil2B\Phi_R/L_R\rceil+1$.  The terms with $N\le N_R$
contribute at most $N_R/\Phi_R\le2B/L_R+2/\Phi_R$.  For the remaining
terms, $x^N\widetilde G_N\le M_0\rho_R^N$ and $\widetilde
Z\ge\ee^{-B\Phi_R}$, so they contribute at most
\begin{align*}
 &\Phi_R^{-1}M_0\ee^{B\Phi_R}\sum_{N>N_R}N\rho_R^N\\
 &\qquad\le
 4M_0(N_R+1)\rho_R\,\ee^{B\Phi_R-N_RL_R}/\Phi_R\\
 &\qquad\le
 4M_0(N_R+1)\rho_R\ee^{-B\Phi_R}/\Phi_R ,
\end{align*}
which tends to $0$.  This proves the claim.

The functions $U_R$ are convex and locally uniformly bounded, so a
subsequence converges locally uniformly on $\R\times(0,\infty)$ to a convex
function $U$.  By the claim, $U$ does not depend on $q$, so its subgradients
lie in the line $\{0\}\times\R$ and its Monge--Amp\`ere measure
$\operatorname{MA}[U]$, in the sense of Alexandrov
\cite[Chapter~1]{Gutierrez2016}, vanishes.
On the other hand, Monge--Amp\`ere measures are weakly continuous under
local uniform convergence of convex functions
\cite[Lemma~1.2.3]{Gutierrez2016}, so \eqref{eq:scale-invariant} gives
$\operatorname{MA}[U]=C\ee^q\,\dd q\,\dd u\ne0$.  This contradiction proves
the theorem.
\end{proof}

\begin{proof}[Proof of \cref{thm:rank-two-intro}]
If $Q$ contains a line, \cref{lem:pointed} gives~(ii).  Such a $Q$ has no
exposed points, so (i) fails.  Otherwise, by \cref{lem:pointed} and
\cref{thm:no-diffuse}, $\mu$ charges the exposed point at the origin.
\Cref{thm:boundary-intro}\textup{(c)} then gives linearly independent
atoms $g_1,g_2$ with
\[
 \supp\mu\subset\Nzero g_1+\Nzero g_2,
\]
and the linear map sending $g_i$ to $e_i$ gives~(i).
\end{proof}

\begin{remark}[Examples and Open Problems]\label{rem:rank-two-open}
The law of two independent Poisson variables $P_1,P_2$, and more generally
that of $(P_1+P_2,P_2)$, whose support is $\Nzero(1,0)+\Nzero(1,1)$, fall
under~(i); the product of a Gaussian and a Poisson law falls under~(ii).  In the normalization of (i) it remains to
show that every positive series $Z(x,y)=\sum_{m,n\ge0}a_{mn}x^my^n$ with
$a_{00},a_{10},a_{01}>0$ and
\[
 \det D^2_{(\log x,\log y)}\log Z=Cxy
\]
is $a_{00}\exp(\lambda_1x+\lambda_2y)$.  In case~(ii) the expected answer is
the product of a Gaussian and a Poisson family.  Neither statement is proved here in general;
\cref{sec:harmonic-boundary} treats three regular regimes.
\end{remark}

\subsection{The Axisymmetric Harmonic Reduction}
\label{sec:harmonic-boundary}

A partial Legendre
transform in the tangential variable turns \eqref{eq:E} into the Laplace
equation for an axisymmetric function on $\R^3$, and positivity of $\mu$
then restricts the boundary behaviour of that function.

From here on, $d=2$ and $b\ne0$.  By an invertible affine change
of the statistic, including a shear of the tangential coordinate, we write
$\theta=(t,s)$ and normalize
\begin{equation}\label{eq:hb-normal}
 \det K''(t,s)=C\ee^t,
 \qquad C>0.
\end{equation}
The partial Legendre transform in the tangential variable is
\begin{equation}\label{eq:hb-legendre}
 p=K_s(t,s),\qquad w(t,p)=K(t,s)-ps.
\end{equation}
Strict convexity makes $(t,s)\mapsto(t,p)$ a local diffeomorphism, and
implicit differentiation gives
\[
 w_p=-s,\qquad w_{pp}=-K_{ss}^{-1},\qquad
 w_{tt}=K_{tt}-\frac{K_{ts}^2}{K_{ss}}=\frac{C\ee^t}{K_{ss}}.
\]
Thus $w_{tt}+C\ee^tw_{pp}=0$.  With
\begin{equation}\label{eq:hb-rv}
 r=2\sqrt C\,\ee^{t/2},\qquad v(r,p)=-w_{pp}=K_{ss}^{-1},
\end{equation}
this becomes
\begin{equation}\label{eq:hb-harmonic}
 w_{rr}+r^{-1}w_r+w_{pp}=0,\qquad
 v_{rr}+r^{-1}v_r+v_{pp}=0.
\end{equation}
The second equation in \eqref{eq:hb-harmonic} says that $v$, the reciprocal
of the tangential variance $K_{ss}$, is harmonic as a function on $\R^3$
that depends only on the distance $r$ to an axis and on the coordinate $p$
along it; such a function is called \emph{axisymmetric}.  We denote by
$\Omega\subset(0,\infty)\times\R$ the image of $(t,s)\mapsto(r,p)$.

The reduction goes back to the Gibbons--Hawking
ansatz~\cite{GibbonsHawking1978}; see Donaldson
\cite[Section~3]{Donaldson2009GGT} and compare his treatment of toric surfaces
\cite[Section~6.1]{Donaldson2009}.  For \eqref{eq:E} in the plane, the relation
to the Laplace equation in $\R^3$ is also noted by Song, Wu and
Xu~\cite{SongWuXu2026}.  What we use in addition is that $K$ is the cumulant
function of a positive measure.  With it we establish rigidity in three
regular boundary regimes.

\begin{theorem}[Regular Boundary Regimes]\label{thm:harmonic-intro}
Let $d=2$ and let $\mu$ be a non-Gaussian solution of \eqref{eq:E}.
\begin{enumerate}
 \item[\textup{(a)}] If $\Theta=\R^2$ and the reciprocal tangential variance
 is the potential of finitely many positive point charges on the axis, of
 arbitrary sizes, then the family is Poisson--Gaussian
 (\cref{thm:finite-center}).
 \item[\textup{(b)}] In the lattice-orthant branch, let
 $F(u)=\sum_{j\ge0}\mu_0\{j\}u^j$ be the generating function of the bottom
 row $\mu_0$, with mean $m=uF'/F$ and variance $V=um'$.  Suppose that $F$
 and the generating function of the adjacent row $\mu_1$ are entire, that
 $m$ maps $(0,\infty)$ onto itself, and that $1/V$ is a Stieltjes function
 of the mean,
 \[
  \frac1V=\frac1m+c+\int_{(0,\infty)}\frac{\sigma(\dd q)}{m+q}
 \]
 with $c\ge0$ and a positive measure $\sigma$, as in
 \cref{thm:stieltjes}.  Then the bottom row is Poisson.
 \item[\textup{(c)}] In the lineality branch, let $\Omega$ be the image of
 $\Theta$ under the partial Legendre map $\theta\mapsto(r,p)$.  No part of
 the boundary of $\Omega$ is a $C^2$ curve $r=R(p)$, $0<R(p)<\infty$, with
 $\Omega$ lying locally on the side $r<R(p)$ (\cref{thm:wall}).
\end{enumerate}
\end{theorem}

The next three subsections prove the three parts of
\cref{thm:harmonic-intro}, which concern three kinds of boundary
behaviour of $v$.  In
\cref{sec:finite-axial}, $v$ is the potential of finitely many point
charges on the axis.  We show that the family then has a limit at the axis,
a one-dimensional NEF whose mean is a meromorphic function of the natural
parameter, and apply the results of Bar-Lev, Bshouty and Enis on such
families \cite{BarLevBshoutyEnis1991}.  \Cref{sec:stieltjes-rows} concerns the
lattice-orthant branch.  There the reciprocal variance of the bottom row,
which is the boundary value of $v$ on the axis, is allowed to be the
potential of a positive measure of charges; it need not be rational, and the
critical-point theorem of Bergweiler and Eremenko
\cite{BergweilerEremenko1995} takes the place of those results.
\Cref{sec:walls} concerns the lineality branch.  There, if $v$
extends to a positive harmonic function on all of $\R^3$, then it is constant
by Liouville's theorem and $K$ is the cumulant function
\eqref{eq:hb-gaussian-poisson} of the product of a Gaussian and a Poisson
family; we show that $\Omega$ cannot end at a finite smooth outer wall.

\subsection{Finite Axial Data}
\label{sec:finite-axial}

The simplest positive axisymmetric harmonic functions are the potentials
\eqref{eq:hb-finite-centers} of finitely many point charges on the axis,
possibly plus a constant; these are the Gibbons--Hawking potentials with
collinear centres.  We refer to such $v$ as \emph{finite axial data}.  For
them the classification is complete.

\begin{theorem}[Finite Axial Data]\label{thm:finite-center}
Suppose that $\Theta=\R^2$ and that
\begin{equation}\label{eq:hb-finite-centers}
 v(r,p)=c+\sum_{j=1}^M
 \frac{q_j}{\sqrt{r^2+(p-a_j)^2}},\qquad
 c\ge0,\quad q_j>0,\quad a_1<\cdots<a_M.
\end{equation}
Then $M\le1$ and the family is Poisson--Gaussian.  After an invertible
affine change of the statistic, it is the product of a Gaussian and a
Poisson family if $M=0$, and of two Poisson families if $M=1$.
\end{theorem}

Since $s=-w_p$ and $s_p=v$, integration gives
\begin{equation}\label{eq:hb-axis-coordinate}
 s=cp+A\log r+B+
 \sum_{j=1}^Mq_j\operatorname{arsinh}\frac{p-a_j}{r}.
\end{equation}
Here the radial integration term is $A\log r+B$ by
\eqref{eq:hb-harmonic}.  A shear of $(t,s)$ changes $A$.  For any component
$I$ of $\R\setminus\{a_1,\ldots,a_M\}$ the shear can be chosen so that $s$
has a finite limit as $r\to0$ with $p\in I$ fixed; we then say that $I$ has
been \emph{selected}.  On a selected component, letting $r\to0$ gives the
derivative of the tangential natural parameter with respect to the
tangential mean:
\begin{equation}\label{eq:hb-axis-vf}
 \frac{\dd s}{\dd p}=v_0(p)
 =c+\sum_{a_j<p}\frac{q_j}{p-a_j}
   +\sum_{a_j>p}\frac{q_j}{a_j-p}.
\end{equation}

This limit comes from an exponential family on the axis, which we call the
\emph{boundary NEF} on $I$.

\begin{lemma}[Axis Limit as a Boundary NEF]\label{lem:hb-axis-limit}
After selecting $I$ as above, \eqref{eq:hb-axis-vf} is the reciprocal of
the variance function of a one-dimensional NEF with mean parameter $p$ and
mean domain $I$.
\end{lemma}

\begin{proof}
Choose the shear so that the coefficient of $\log r$ in
\eqref{eq:hb-axis-coordinate} cancels on $I$.  Then $p(t,s)$ converges
locally uniformly as $t\to-\infty$ to the inverse of the strictly
increasing limiting map $s=s(p)$.  Fix $s_0$.  For small $h$, the ratios
\[
 \frac{Z(t,s_0+h)}{Z(t,s_0)}
\]
are moment-generating functions of the tangential marginal under the
$(t,s_0)$-tilted law.  Local convergence of their logarithmic derivatives,
together with bounds at $h=\pm\varepsilon$, gives tightness and uniform
exponential integrability.  Their limit is therefore a moment-generating
function.  Its first two logarithmic derivatives are $p$ and $1/v_0(p)$.
Because $\Theta=\R^2$, every sheared ray in this construction remains in the natural
domain, and the mean domain is the whole selected component.
\end{proof}

The next lemma shows that a one-dimensional NEF is supported on a lattice
if, near a finite endpoint of its mean domain, the natural parameter is a
multiple of the logarithm of the distance from the mean to the endpoint plus
an analytic function of the mean.

\begin{lemma}[Endpoint Lattice Coordinate]\label{lem:hb-endpoint-lattice}
Let $Z_0$ be the Laplace transform of a nonzero positive measure supported
in $[a,\infty)$ and finite for every real natural parameter.  Suppose that
the natural parameter $s$ and the mean $m$ are related by
\[
 s=q\log(m-a)+h(m),\qquad q>0,
\]
near $a$, where $h$ is real analytic at $a$.  Then, for
$u=\ee^{s/q}$,
\begin{equation}\label{eq:hb-endpoint-series}
 \ee^{-as}Z_0(s)=F(u),
\end{equation}
where $F$ is entire with nonnegative Taylor coefficients.  Equivalently,
the measure is supported on $a+q^{-1}\Nzero$.
\end{lemma}

\begin{proof}
Analytic inversion gives $m-a=\Phi(u)$ with $\Phi(0)=0$ and
$\Phi'(0)>0$.  Since $\dd s=q\,\dd u/u$,
\[
 \log\bigl(\ee^{-as}Z_0(s)\bigr)
 =q\int_0^u\frac{\Phi(z)}z\,\dd z+\text{constant}
\]
is analytic at zero.  Push the measure forward by
$x\mapsto q(x-a)$ and write $u=\ee^{-\tau}$.  Uniqueness of the Laplace
transform in $\tau$ identifies this pushforward with the discrete measure
given by the Taylor coefficients.  They are therefore nonnegative.  Because
$Z_0$ is finite for every real $s$, the series is finite for every
positive $u$ and hence is entire.
\end{proof}

\begin{proof}[Proof of \cref{thm:finite-center}]
If $M\ge1$, apply \cref{lem:hb-axis-limit} on the right
exterior component $(a_M,\infty)$.  There
\[
 \frac1{V(p)}=v_0(p)=c+\sum_{j=1}^M\frac{q_j}{p-a_j}.
\]
Near the endpoint, $s=q_M\log(p-a_M)+h(p)$ with $h$ analytic.  With
$u=\ee^{s/q_M}$, the boundary NEF's mean domain $(a_M,\infty)$ makes its
closed convex support $[a_M,\infty)$.  \Cref{lem:hb-endpoint-lattice}
therefore shows that its transform is $\ee^{a_Ms}F(u)$, where $F$ is entire
with nonnegative coefficients.  Its mean
$p(s)=a_M+uF'(u)/(q_MF(u))$ is therefore a nonconstant meromorphic function
of the complexified natural parameter $s\in\mathbb C$, and it satisfies
$p'=V(p)$, first for real $s$ and then on $\mathbb C$ by the identity
theorem.  Bar-Lev, Bshouty and Enis proved that if $W$ is meromorphic on
$\mathbb C$ and not identically zero, and the equation $m'=W(m)$ has a
nonconstant meromorphic solution on $\mathbb C$, then $W$ is a polynomial of
degree at most two \cite[Theorem~4.1]{BarLevBshoutyEnis1991}.  Since
$V=1/v_0$ is rational, $V$ is therefore a polynomial of degree at most two.
If $c>0$, then
$V(p)\to c^{-1}$ as
$p\to\infty$, while $V$ has a simple zero at $a_M$, a contradiction.
Thus $c=0$.  Now $V(p)\sim p/(\sum_jq_j)$ at infinity, so $V$ is linear;
comparison of reciprocals with the partial fractions above leaves only the
single pole $M=1$.

For the surviving single pole, integrate \eqref{eq:hb-axis-coordinate} and
invert the partial Legendre transform.  After an affine change of $(t,s)$
that preserves \eqref{eq:hb-normal} and the addition of an affine function
to $K$,
\begin{equation}\label{eq:hb-two-poisson}
 K(t,s)=A\ee^{\alpha t+\beta s}
       +B\ee^{(1-\alpha)t-\beta s},\qquad
 A,B>0,\quad\beta\ne0.
\end{equation}
The two exponent vectors are
independent, so \eqref{eq:hb-two-poisson} is a product of two Poisson
families after an invertible linear change of the statistic.

If $M=0$, positivity gives $c>0$.  From $w_{pp}=-c$ and
\eqref{eq:hb-harmonic}, followed by a shear and the removal of affine
terms,
\[
 w=-\frac c2p^2+cC\ee^t.
\]
Legendre inversion yields
\begin{equation}\label{eq:hb-gaussian-poisson}
 K(t,s)=\frac{(s+\beta)^2}{2c}+cC\ee^t+\gamma,
\end{equation}
the cumulant function of the product of a Gaussian and a Poisson family.
This proves both alternatives.
\end{proof}

\subsection{Regular Stieltjes Rows}
\label{sec:stieltjes-rows}

In the lattice-orthant branch of \cref{thm:rank-two-intro}, use the
coordinates of part~\textup{(i)} of that theorem, in which
$\supp\mu\subset\Nzero^2$ and $b=(1,1)$.  Let
$F(u)=\sum_{j\ge0}\mu_0\{j\}u^j$, with $u=\ee^s$, be the generating function
of the bottom row, so that $G_0(s)=F(\ee^s)$; let $m=uF'/F$ be the mean of
the NEF generated by $\mu_0$, and write $V=um'$ for its variance.  As
$t\to-\infty$ the laws $P_{(t,s)}$ concentrate on the bottom row and the
tangential variance tends to $V$, so $1/V$ is the boundary value on the
axis $r=0$ of the harmonic function $v$ of \eqref{eq:hb-rv}.  Taking Laplace
transforms in \eqref{eq:leading-row} gives $F^2VG_1=CuF^3$; hence the
Laplace transform of the first row, as a function of $u$, is
\begin{equation}\label{eq:hb-first-row}
 G_1=C\frac{uF(u)}{V(m(u))}=C\frac{F(u)}{m'(u)}.
\end{equation}
The hypothesis in \cref{thm:stieltjes} that $F/m'$ is entire
therefore says exactly that the Laplace transform of the adjacent row
$\mu_1$ is finite for every tangential parameter.  This is an additional
hypothesis, not a consequence of \cref{thm:rank-two-intro}.  The
hypothesis on the variance is that $1/V-1/m$ is a nonnegative constant plus
the Stieltjes transform of a positive measure $\sigma$, so that $1/V$ is a
Stieltjes function of the mean
\cite[Chapter~2]{SchillingSongVondracek2012}.  For finitely supported
$\sigma$ this is the form \eqref{eq:hb-axis-vf} taken by finite axial data
on the component $(a_M,\infty)$, with $a_M=0$ and $q_M=1$.

\begin{theorem}[Stieltjes Rows]\label{thm:stieltjes}
Let $F(u)=\sum_{n\ge0}A_nu^n$ be entire, where $A_n\ge0$, $F(0)=1$ and
$A_1>0$.  Put $m(u)=uF'(u)/F(u)$ and $V(m(u))=um'(u)$.  Suppose that
$m:(0,\infty)\to(0,\infty)$ is onto and that, for a positive Borel measure
$\sigma$ on $(0,\infty)$,
\[
 \frac1{V(m)}=\frac1m+h(m),\qquad
 h(m)=c+\int_{(0,\infty)}\frac{\sigma(\dd q)}{m+q},
\]
where $c\ge0$ and $c+\int_{(0,\infty)}\sigma(\dd q)/q<\infty$.
If $F/m'$ extends to an entire function, then $c=0$, $\sigma=0$, and
$F(u)=\exp(A_1u)$.
\end{theorem}

With $h$ as in the theorem, set
\begin{equation}\label{eq:hb-stieltjes-data}
 \mathcal H(m)=\int_0^m h(x)\,\dd x,
 \qquad J(m)=m+\int_0^m xh(x)\,\dd x.
\end{equation}
These functions parametrize $u$ and $F$ by the mean $m$ along the positive
axis.

\begin{lemma}[Stieltjes Parametrization]\label{lem:hb-stieltjes-param}
Under the hypotheses of \cref{thm:stieltjes},
\begin{equation}\label{eq:hb-stieltjes-param}
 A_1u=m\ee^{\mathcal H(m)},\qquad F(u(m))=\ee^{J(m)}
\end{equation}
on the positive axis.
\end{lemma}

\begin{proof}
Since $V=\dd m/\dd\log u$,
\[
 \frac{\dd\log u}{\dd m}=\frac1m+h(m).
\]
The first identity follows by integration and
$m=A_1u+O(u^2)$.  Moreover
\begin{equation*}
 \dd\log F/\dd m=m\,\dd\log u/\dd m=1+mh(m),
\end{equation*}
which gives the second.
\end{proof}

The parametrization controls the order of the entire function $F$.

\begin{lemma}[Order Less than One]\label{lem:hb-subunit-order}
If $h\not\equiv0$, then the order of the entire function $F$ is strictly
smaller than one.  More precisely, with $M_\sigma=\sigma((0,\infty))$, the order is zero
if $c>0$ or $M_\sigma=\infty$, and is at most $(1+M_\sigma)^{-1}$ if
$c=0$ and $0<M_\sigma<\infty$.
\end{lemma}

\begin{proof}
Nonnegative coefficients give $\max_{|u|=R}|F(u)|=F(R)$, so it suffices to
use \eqref{eq:hb-stieltjes-param} on the positive axis.  If $c>0$, then
$\mathcal H(m)\ge cm$, whence $m=O(\log u)$; since
$h(m)\le h(0)<\infty$, one has
$J(m)=O(m^2)$ and the order is zero.

Now suppose $c=0$ and put $M_T=\sigma((0,T])$.  Tonelli's theorem gives
\begin{equation}\label{eq:hb-H-lower}
 \mathcal H(m)=\int\log\left(1+\frac mq\right)\sigma(\dd q)
 \ge M_T\log(1+m/T),
\end{equation}
so $m=O_T(u^{1/(1+M_T)})$.  If $M_\sigma<\infty$, then
$J'(m)\le1+M_\sigma$; letting $M_T\uparrow M_\sigma$ gives order at most
$(1+M_\sigma)^{-1}<1$.  If $M_\sigma=\infty$, the integrability condition
$\int\sigma(\dd q)/q<\infty$ of \cref{thm:stieltjes}
and $0\le y-\log(1+y)\le y^2/2$ give
\[
 J(m)=m+\int\left[m-q\log(1+m/q)\right]\sigma(\dd q)=O(m^2).
\]
Because $M_T$ can be arbitrarily large, the order is at most
$2/(1+M_T)$ for every $T$, and is therefore zero.
\end{proof}

\begin{proof}[Proof of \cref{thm:stieltjes}]
The function $m(u)=uF'(u)/F(u)$ is meromorphic.  Since $F/m'$ is entire,
$m'$ has no finite zero: this is immediate where $F\ne0$; at a nonzero
zero $u_0$ of $F$ of multiplicity $k$,
\[
 m(u)=\frac{ku_0}{u-u_0}+O(1),\qquad
 m'(u)=-\frac{ku_0}{(u-u_0)^2}+O(1),
\]
and $m'(0)=A_1>0$.

Suppose $h\not\equiv0$.  The order of a logarithmic derivative does not
exceed that of the function, so \cref{lem:hb-subunit-order} gives
$\operatorname{ord}(m)\le\operatorname{ord}(F)<1$.  The derivative of a
transcendental meromorphic function of order less than one has infinitely
many zeros~\cite[Theorem~4(a)]{BergweilerEremenko1995}; hence $m$ is
rational.  Every finite pole of $m$ is simple.  If a rational map of degree
$k\ge2$ had no finite critical point, the point at infinity would have to
account for all $2k-2$ critical points counted by the Riemann--Hurwitz
formula, although it contributes at most $k-1$.  Thus $m$ is a M\"obius
transformation.  It is real and positive on the positive
axis, $m(0)=0$, $m'(0)=A_1$, and it maps $(0,\infty)$ onto itself.  Therefore
$m(u)=A_1u$: the alternative $A_1u/(1+bu)$ is bounded for $b>0$ and has a
positive-axis pole for $b<0$.  Consequently $F=\exp(A_1u)$ and $V=m$.
The representation of $1/V$ assumed in the theorem then forces $c=0$ and
$\sigma=0$, contradicting $h\not\equiv0$.  Hence $h\equiv0$, that is, $c=0$
and $\sigma=0$, and \eqref{eq:hb-stieltjes-param} then gives $A_1u=m$ and
$F(u)=\ee^{m}=\exp(A_1u)$.
\end{proof}

\subsection{Walls in the Lineality Branch}
\label{sec:walls}

In the lineality branch the convex support $Q$ is invariant under tangential
translations.  We show that the image $\Omega$ of the partial Legendre map
then cannot end, on its outer side, at a smooth curve.  We call a graph
\begin{equation}\label{eq:hb-wall}
 \mathcal W=\{(R(p),p):p\in I\},\qquad 0<R(p)<\infty,
\end{equation}
a \emph{finite $C^2$ outer radial wall} if $R\in C^2(I)$ and $\Omega$ lies
locally
on the side $r<R(p)$, with $\mathcal W\subset\partial\Omega$, and, after
shrinking $I$ if necessary, there is an $\varepsilon>0$ such that
$\varepsilon<\inf_I R$ and
\[
 \{(r,p):p\in I,\ R(p)-\varepsilon<r<R(p)\}\subset\Omega.
\]

A positive harmonic function is integrable up to a smooth piece of the
boundary.

\begin{lemma}[Wall Integrability]\label{lem:hb-wall-integrability}
For every $I'\Subset I$ and $0<\delta<\varepsilon$,
\[
 \int_{I'}\int_{R(p)-\delta}^{R(p)}v(r,p)\,\dd r\,\dd p<\infty.
\]
\end{lemma}

\begin{proof}
Lift $v$ to the positive harmonic function
$v(\sqrt{X^2+Y^2},p)$ in the corresponding three-dimensional domain.
Choose a bounded $C^2$ subdomain whose boundary contains a compact subpatch
of the surface of revolution generated by $\mathcal W$.  A positive harmonic
function there is the Poisson integral of a finite positive boundary
measure.  The integral of the Poisson kernel over the domain is uniformly
bounded in the boundary point; equivalently, it is the inward normal
derivative of the Dirichlet solution of $-\Delta q=1$.  Tonelli's theorem
gives local $L^1$ integrability in three dimensions.  On $r\ge r_0$,
cylindrical volume is equivalent to $\dd r\,\dd p\,\dd\vartheta$; apply
this with $r_0=\inf_I R-\varepsilon>0$ to prove the assertion.
\end{proof}

Put $n=K_t$, the mean of the normal coordinate.  At every positive radius,
\begin{equation}\label{eq:hb-wall-identities}
 w_r=\frac{2n}{r},\qquad n_r=\frac r2v.
\end{equation}
Indeed $n=w_t=(r/2)w_r$, and the second identity follows by differentiating
and using \eqref{eq:hb-harmonic}.

By these identities, a point of a wall would give a boundary point of the
mean domain that is approached with bounded relative entropy.
\Cref{thm:massive-face} would then produce an exposed face, which
lineality excludes.

\begin{theorem}[No Finite Smooth Wall]\label{thm:wall}
If $Q$ contains a line in the tangential direction, then $\Omega$ has no
finite $C^2$ outer radial wall.
\end{theorem}

\begin{proof}
The set $Q$ is invariant under all translations along the line.  If a wall
\eqref{eq:hb-wall} existed, choose
$p_0\in I$ for which \cref{lem:hb-wall-integrability} gives an
integrable radial slice.  By \eqref{eq:hb-wall-identities}, as
$r\uparrow R(p_0)$ the normal mean $n$ and then $w$ have finite limits.
Also $t=2\log(r/(2\sqrt C))$ remains finite, and
\begin{equation}\label{eq:hb-dual-wall}
 K^*(n,p_0)=tn-w
\end{equation}
has a finite limit.  Thus the means $m=(n,p_0)$ and the values $K^*(m)$ of the dual potential
of \cref{sec:boundary} converge to finite limits.

If $\liminf|s|<\infty$, take a sequence with $s$ bounded.  Then $(t,s)$
converges finitely and $K=w+p_0s$ stays finite.  Fatou's lemma gives a
finite limiting Laplace transform, and \cref{thm:regularity-intro}
puts the limit in $\Theta$.  But the Jacobian of
$(t,s)\mapsto(r,p)$ is $(r/2)K_{ss}>0$, so $(R(p_0),p_0)$ is an interior
point of $\Omega$, a contradiction.

Otherwise $|s|\to\infty$; pass to a subsequence of fixed sign.  The finite
mean limit lies in $\partial M$ (otherwise $\nabla K^*$ keeps $(t,s)$
finite), and the normalized natural parameters tend to $(0,1)$ or
$(0,-1)$.  \Cref{thm:massive-face}, applied using
\eqref{eq:hb-dual-wall}, produces a proper exposed face with this tangential
normal.  Tangential lineality makes the support function infinite in both
nonzero tangential directions, so no such face exists.  This final
contradiction excludes the wall.
\end{proof}

\Cref{thm:wall} leaves open the other ways in which $\Omega$ can end
in the lineality branch: at the symmetry axis $r=0$, at infinite radius, at
a singular or cuspidal wall, or with the tangential mean leaving every
compact interval.

\printbibliography[heading=bibintoc]

\end{document}